\documentclass[]{article}
\usepackage[noadjust]{cite}
\usepackage{latexsym,amssymb,enumerate,amsmath,epsfig,amsthm,authblk,dsfont,fullpage}
\usepackage{multirow}
\usepackage{setspace,color,soul}
\usepackage{tikz}
\usepackage{algorithm,algorithmic}
\usetikzlibrary{shapes,arrows}
\usepackage[normalem]{ulem}
\usepackage{textcomp}
\usepackage{tikz-cd}
\usepackage{accents}
\usepackage{appendix}
\usepackage{graphicx,subfigure}
\usepackage{makecell}
\usepackage{tablefootnote} 
\usepackage{epstopdf}
\usepackage{grffile}
\usetikzlibrary{decorations.pathreplacing}
\usepackage{wrapfig}
\usepackage{hyperref}
\hypersetup{
	colorlinks=true,
	linkcolor=blue,
	citecolor=blue,
	filecolor=blue,      
	urlcolor=blue,
}

\newcommand{\commentout}[1]{}

\newcommand{\resid}{\text{resid}}
\newcommand{\bc}{\mathbf{c}}

\newcommand{\wT}{\widetilde{T}}
\newcommand{\by}{\mathbf{y}}
\newcommand{\bx}{\mathbf{x}}

\newcommand{\bA}{\mathbf{A}}

\newcommand{\bF}{\mathbf{F}}
\newcommand{\bff}{\mathbf{f}}

\newlength{\dhatheight}

\begin{document}
\title{
Group projected Subspace Pursuit for Identification of variable coefficient differential equations (GP-IDENT)}

\author{Yuchen He \thanks{Institute of Natural Sciences, Shanghai Jiao Tong University. Email: yuchenroy@sjtu.edu.cn (Yuchen He is the corresponding Author.  The author list is alphabetized.)}, Sung-Ha Kang \thanks{School of Mathematics, Georgia Institute of Technology. Email: kang@math.gatech.edu. Research is supported in part by Simons Foundation 584960.}, Wenjing Liao \thanks{School of Mathematics, Georgia Institute of Technology. Email: wliao60@gatech.edu. Research is supported in part by NSF grant NSF-DMS 2145167.}, 
Hao Liu \thanks{Department of Mathematics, Hong Kong Baptist University. Email: haoliu@ hkbu.edu.hk. Research is supported in part by HKBU 162784, HKBU 179356 and NSFC 12201530.}
and
Yingjie Liu \thanks{School of Mathematics, Georgia Institute of Technology. Email: yingjie@math.gatech.edu. }}

\date{}

\maketitle
\begin{abstract}
We propose an effective and robust algorithm for identifying partial differential equations (PDEs) with space-time varying coefficients from a single  trajectory of noisy observations.  Identifying unknown differential equations from noisy observations is a difficult task, and it is even more challenging with space and time varying coefficients in the PDE.    The proposed algorithm, GP-IDENT, has three ingredients: (i) we use B-spline bases to express the unknown space and time varying coefficients, (ii) we propose Group Projected Subspace Pursuit (GPSP) to find a sequence of candidate PDEs with various levels of complexity, and (iii) we propose a new criterion for model selection using the Reduction in Residual (RR) to choose an optimal one among the pool of candidates.   The new GPSP considers group projected subspaces which is more robust than existing methods in distinguishing correlated group features. 
We test GP-IDENT on a variety of PDEs and PDE systems, and compare it with the state-of-the-art parametric PDE identification algorithms under different settings to illustrate its outstanding performance.
Our experiments show that GP-IDENT is effective in identifying the correct terms from a large dictionary and the model selection scheme is robust to noise. 
	\end{abstract}

\section{Introduction}
Partial Differential Equations (PDEs) are indispensable and ubiquitous mathematical method articulating fundamental laws that govern various phenomena in physics, chemistry, biology,  finance, and many other fields.   
Let the variable of the given data be
$u(x,t): \Omega\times[0,T_{\max}]\to\mathbb{R}$, where $\Omega\subset\mathbb{R}^{d}$ is a $d$-dimensional  spacial domain, and $T_{\max}>0$ is the final time of the observation. An important class of models that describe the dynamical features of $u$ is the evolution PDE~\cite{galaktionov2003stability,carvalho2012attractors,rudnicki2004chaos}
	\begin{align}
		u_t = \mathcal{F}(u,\partial_xu, \partial_x^2u,\cdots)\label{eq_PDE_general}
	\end{align}
with a functional $\mathcal{F}$. In the multidimensional case with $d>1$, the spatial location is given by $x=(x_1,\dots,x_{d})$.  We use the multi-indexing notation $\partial_x^mu=\{\partial_x^\alpha u:=\partial^{\alpha_1}_{x_1}\partial^{\alpha_2}_{x_2}\cdots\partial^{\alpha_{d}}_{x_{d}}u\;,\alpha = (\alpha_1,\dots,\alpha_{d})\;, \alpha_{1}+\cdots+\alpha_{d}=m\}$ to denote the collection of $m$-order partial derivatives of $u$. 
 The model~\eqref{eq_PDE_general} covers a wide range of important PDEs including the advection-diffusion equation for transferring physical quantities, the Kolmogorov-Petrovsky-Piskunov (KPP) equation~\cite{tikhomirov1991study} for population genetics, the incompressible Navier-Stokes equation~\cite{majda2002vorticity}, the Korteweg-de Vries (KdV) equation~\cite{newell1985solitons}, and the Kuramoto-Sivashinsky (KS) equation~\cite{kuramoto1978diffusion} for fluid dynamics. For a vector-valued $u$, \eqref{eq_PDE_general} also covers PDE systems such as the nonlinear Schr\"{o}dinger equation~\cite{zakharov1974complete} for light propagation. 
Model \eqref{eq_PDE_general} can be regarded as an infinite dimensional dynamical system whose asymptotic properties such as attractors~\cite{carvalho2012attractors} and chaotic behaviors~\cite{rudnicki2004chaos} have been extensively studied.

Classical approaches to derive PDE for specific physical processes are based on physical laws and simplified assumptions.  In modern science, 
\textbf{data-driven PDE identification} 
is explored which automatically identifies such model~\eqref{eq_PDE_general} from the given observation. Such approaches allow scientists and engineers to discover non-linear and high-order complicated PDEs which are hard to model by empirical experience.

In literature, various techniques have been developed to identify the active features, where sparse regression is one of the major frameworks for PDE identification~\cite{brunton2016discovering,schaeffer2017learning,rudy2017data,kang2021ident,kaiser2018sparse,loiseau2018constrained,mangan2017model,messenger2021weak,Rudy2019DataDrivenIO,li2020robust,zhang2018robust,tang2022weakident}.
Brunton et al.~\cite{brunton2016discovering} studied the application of $L_1$-norm regularization in the context of PDE identification and proposed the sequential thresholded least-squares  
to find the active features. Kang et al.~\cite{kang2021ident} proposed to obtain a series of candidate models using LASSO~\cite{tibshirani1996regression} and then select the optimal model with the minimal time evolution error (TEE).  
Rudy et al.~\cite{rudy2017data} penalized the coefficients using  the $L_0$-norm, and proposed sequential threshold ridge regression (STRidge) to solve the resulting problem. 
He et al.~\cite{he2022robust} proposed to use Subspace Pursuit (SP)~\cite{dai2009subspace} with a series of sparsity levels to generate candidate models. They also proposed Successively Denoised Differentiation (SDD) for denoising the input, and cross-validation error evaluation 
and multi-shooting TEE for selecting the optimal candidate.
Other sparsity promoting penalties are studied in ~\cite{kaptanoglu2021promoting,champion2020unified,carderera2021cindy}. A theoretical analysis for PDE identification can be found in \cite{he2022much,zhang2019convergence, he2022asymptotic}. 
Methods such as~\cite{messenger2021weakPDE,messenger2021weak,tang2022weakident} used sequential least squares \cite{messenger2021weakPDE,messenger2021weak} and subspace pursuit \cite{tang2022weakident} for a weak form of PDE instead of differential form which are more robust to noise. 

Another line of works are based on neural networks~\cite{long2018pde,wu2020data,xu2021deep,rao2022discovering}, where sparse regression is embedded for feature selection, and a sufficient amount of trajectories of data are required for training. Different frameworks such as symbolic regression~\cite{bongard2007automated,schmidt2009distilling,vaddireddy2019equation,long2019pde} are also available.  See \cite{north2022review} for a recent review.

PDEs with space and time  varying coefficients are often used in real applications, such as optimal control~\cite{vazquez2008control,izadi2015pde,kerschbaum2019backstepping}, trajectory planning~\cite{meurer2009trajectory}, studies of piezoelectricity~\cite{kaltenbacher2007identification}, and electromagnetic eddy current problems~\cite{langer2019space}. In such cases,  certain coefficients may depend on both time and space, and in some equations, parts of the coefficients may vary with time while the others vary with space. An effective and robust PDE identification scheme with the flexibility of handling space and time varying coefficients is in need.

\begin{figure}
		\centering
			\includegraphics[width = 0.85\textwidth]{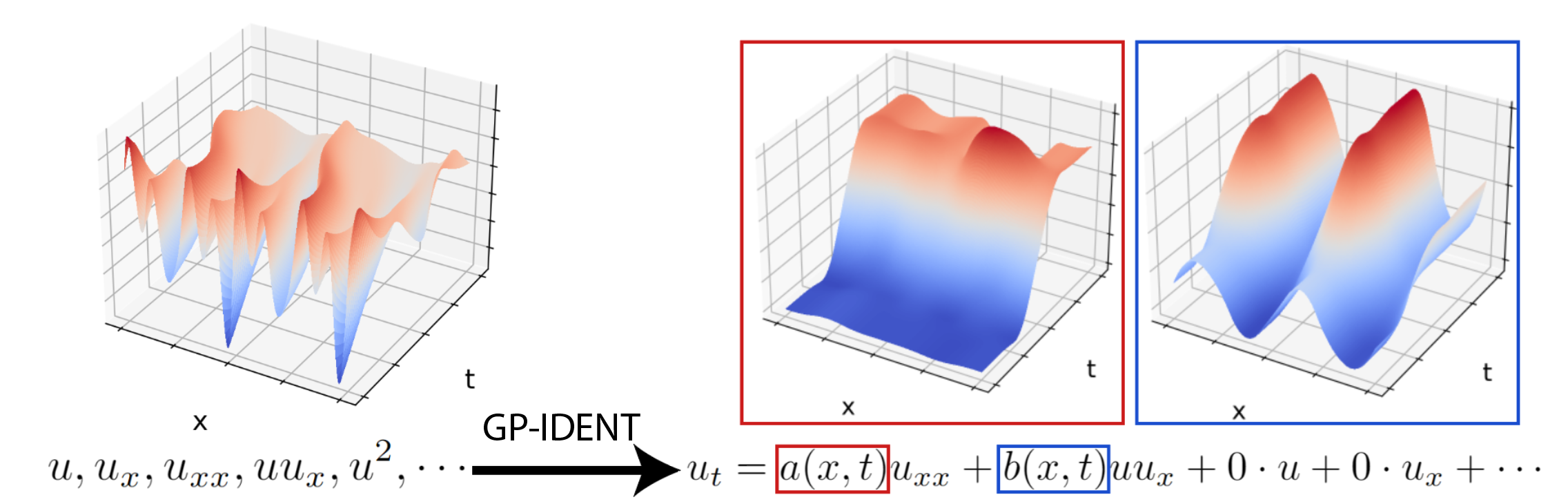}
	\caption{The proposed algorithm, GP-IDENT, identifies the underlying PDE with space-time varying coefficients from a single   trajectory of noisy observations. }\label{fig_GPSP}
\end{figure}

The goal of PDE identification in this paper is to find an expression of~\eqref{eq_PDE_general} in a parametric form
	\begin{align}
	u_t = \sum_{g=1}^GC_g(x,t)f_g(x,t)\label{eq_PDE_spec1}
\end{align}
based on single, possibly noisy, observations of a solution trajectory $u$ in $\Omega\times [0,T_{\max}]$.  The set of potential features $\mathcal{G}=\{f_g\}_{g=1}^G$ forms a dictionary, which can include linear terms such as partial derivatives of $u$ in various orders, and products of multiple linear terms, e.g., $uu_x$ and $u^2$. The size of the dictionary $G>0$ is sufficiently large, and $C_g$, $g=1,2,\dots,G$, represents a  space-time dependent function. 
Figure~\ref{fig_GPSP} provides an illustration: from a noisy observation of a single  solution trajectory, the proposed method identifies the features $u_{xx}$ and $uu_x$ from a dictionary and reconstructs the  respective space and time varying coefficients, i.e., $a(x,t)$ and $b(x,t)$. 

There are  few works dealing with space-time varying coefficients: \cite{kang2021ident,Rudy2019DataDrivenIO} laid out a framework to identify varying coefficients and explored regularizers to encourage structural sparsity.  Algorithms were numerically tested on PDEs with either space or time varying coefficients. In \cite{kang2021ident}, the authors explored identification of spatially varying coefficients with Group-Lasso and proposed a Base Element Expansion (BEE) technique.  
In \cite{he2022numerical}, authors proposed a split Bregman method to identify interacting kernels in aggregation equations, where the kernel to be identified is space and time varying, yet the form of the equation is assumed to be given.

In this paper, we propose Group Projected subspace pursuit for the IDENTification of variable coefficient PDEs (GP-IDENT) to identify parametric PDEs with space-time varying coefficients from a single trajectory of noisy data.
Spanning the hypothesis space by B-spline bases~\cite{schoenberg1988contributions}, our strategy is to generate a collection of candidate models by using different levels of group sparsity, 
then evaluate each candidate by considering the Reduction in Residual error (RR) to identify the optimal model. Since the candidate generation involves solving a non-convex, non-differentiable, NP-hard problem~\cite{dai2009subspace},  we  design a novel and effective Group Projected Subspace Pursuit (GPSP) greedy algorithm to produce candidate models with any specified level of group sparsity. 
We compare these methods on a variety of linear, non-linear PDEs and systems of differential equations with different levels of noise. 
Our experiments show that GP-IDENT outperforms other methods in terms of effectiveness, efficiency, and robustness.

Contribution can be summarized as follows. 
\begin{enumerate}
\item We propose a novel method, GP-IDENT, to identify parametric PDEs with variable coefficients which varies in space and time.  We assume the given data is a single observation of possibly noisy data.  The proposed procedure integrating SDD, GPSP, and RR shows robust performances compared to other state-of-the-art approaches. 
\item We propose a new Group Projected Subspace Pursuit algorithm, GPSP, for structured sparse regression with group $\ell_0$-norm constraint.  
GPSP is efficient in searching for the correct features in the underlying PDE, and outperforms block subspace pursuit 
\cite{kamali2013block} especially when different groups or columns within a group are highly correlated. 
\item We propose to consider the Reduction in the Residual error (RR)  to identify the optimal model, which give more stable identification results compared to AIC-based approaches~\cite{Rudy2019DataDrivenIO,li2020robust} when the data is noisy or the dictionary is large.
\end{enumerate}
This paper is organized as follows. In Section~\ref{sec_prop}, we present the detailed procedure of the proposed method, GP-IDENT.  
In Section~\ref{sec_GPSP}, we describe the new Group Projected Subspace Pursuit algorithm, GPSP, and explain the details, including comparisons with block subspace pursuit~\cite{kamali2013block}.   Following numerical implementation details in Section \ref{Sec:implementation}, in Section~\ref{sec_exp}, we present numerical experiments to validate the effectiveness of the proposed GP-IDENT and compares it with the state-of-the-art methods on various types of PDEs. We conclude the paper with some discussions in Section~\ref{sec_conclusion}.

\section{Group Projected subspace pursuit for IDENTification (GP-IDENT) of variable coefficient differential equations
}\label{sec_prop}

\begin{figure}
		\centering
			\includegraphics[width = 0.75\textwidth]{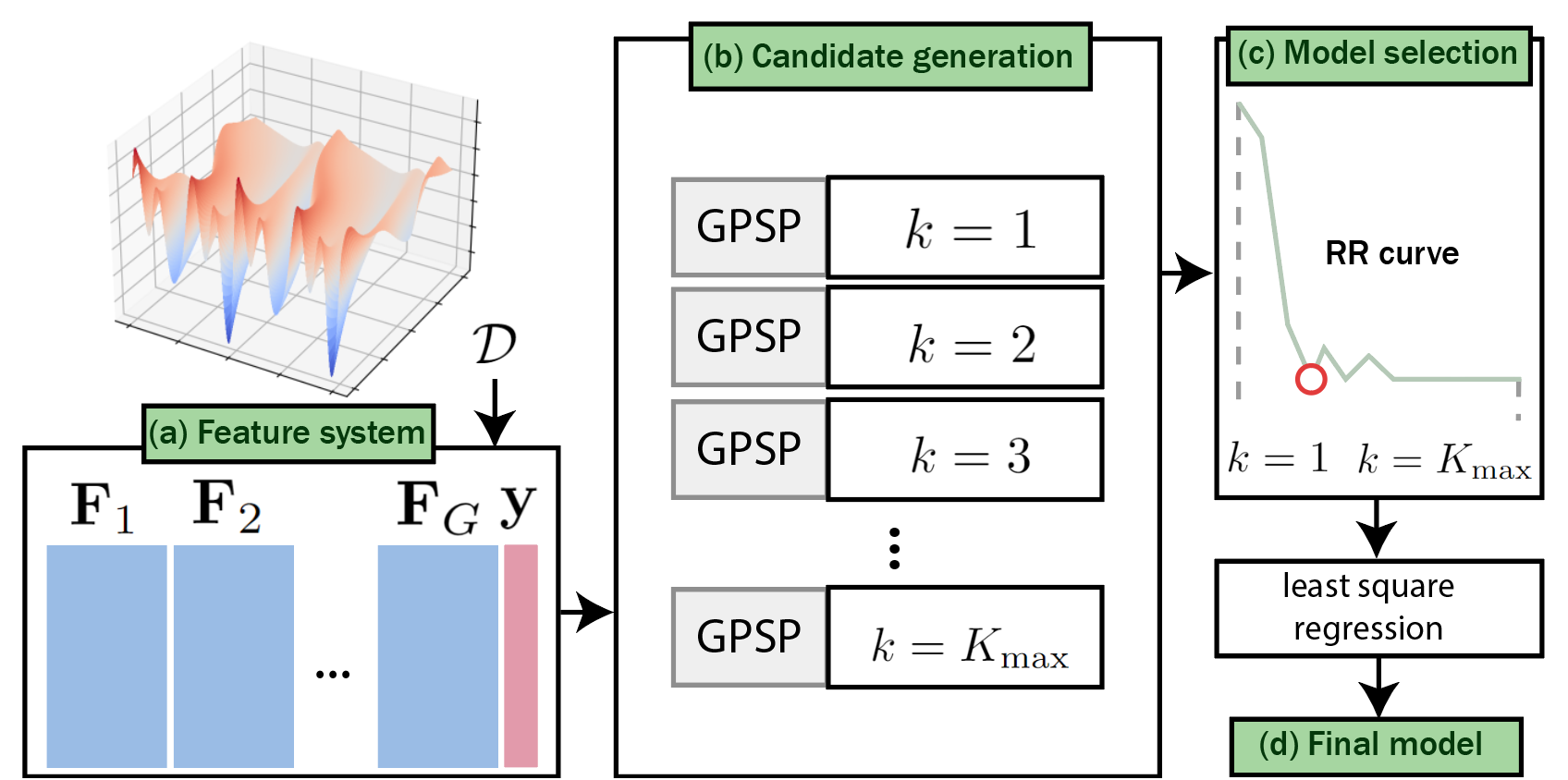}
	\caption{Workflow of the proposed GP-IDENT for varying coefficient PDE identification from noisy observations. (a) Given a noisy observation of a single trajectory, we build the feature system (Section~\ref{ss:feature}). 
 (b) For $k=1,\dots, K_{\max}$, we generate a candidate model 
 by solving a group-$\ell_0$ optimization problem, the proposed GPSP algorithm (Section~\ref{sec_GPSP}). (c) We evaluate each candidate's Reduction in the Residual error (RR) to  select the optimal model among candidate models. 
 (Section~\ref{ss:rss}). (d) 
 Reconstruct the coefficients by least square regression.
}\label{fig_workflow}
\end{figure}

The proposed method has four steps as illustrated by the flowchat in Figure~\ref{fig_workflow}. \textbf{[Step 1]} From the noisy single observation, to account for the instability caused by noise, we employ  the Successively Denoised Differentiation (SDD)~\cite{he2022robust} to smooth the data and generate the feature system as in \cite{Rudy2019DataDrivenIO,kang2021ident}.  Each variable coefficient is represented by B-spline bases~\cite{schoenberg1988contributions} to account for the variation in space and time.  The details are presented in subsection~\ref{ss:feature}. 
\textbf{[Step 2]} To find candidate models of each sparsity level, we propose GPSP. We  describe the procedure in subsection~\ref{ss:candidate} and the algorithmic details are presented in Section~\ref{sec_GPSP}.  
\textbf{[Step 3]} Among the candidate models, we present the model selection criterion based on a Reduction in the Residual error (RR), detailed in subsection~\ref{ss:rss}.  \textbf{[Step 4]} Finally the coefficients are reconstructed.  
We summarize the  proposed GP-IDENT algorithm in 
Algorithm \ref{alg_Prop}. 

\textbf{Notation:}
In this paper, we use standard letters such as $u,B$ for scalars. We use bold lowercase letters such as $\bc$ for vectors, and bold uppercase letters such as $\bA$ for matrices. For a matrix $\bA$, $\bA^\top$ denotes its transpose, and $\bA^\dagger$ denotes its pseudo-inverse. A vector $\bc\in\mathbb{R}^N$ is viewed as a column vector, and $\bc^{\top}$ as its transpose is a row vector. $\|\bc\|_1=\sum_{n=1}^N|c_n|$ and $\|\bc\|_2=\sqrt{\sum_{n=1}^Nc_n^2}$ are $\ell_1$ and $\ell_2$-norm of $\bc$, respectively. We use $\text{supp}(\bc):=\{n=1,2,\dots,N| c_n\neq 0\}$ for the set of indices of the non-zero entries of $\bc$, and its $\ell_0$-norm $\|\bc\|_0$ is the number of elements in  $\text{supp}(\bc)$.

 		\begin{algorithm}[t]
		\begin{algorithmic}[1]
			\REQUIRE{Sampled trajectory data $\mathcal{D}$, over-complete dictionary $\mathcal{G}$, smoothing window size $w\geq 0$, hypothesis space $\mathcal{H}_M$, maximal sparsity level $K_{\max}$, threshold  $\rho$, and selection window $L$
			}
			\STATE Construct the  feature system $(\bA, \by)$ based on $\mathcal{D}$, $\mathcal{H}_M$, and $\mathcal{G}$ using SDD with window size $w$.
			
			\FOR{$k=1,\dots,K_{\max}$}
			\STATE Obtain an 
   approximate solution $\overline{\bc}^*(k)$ with  GPSP (Section~\ref{sec_GPSP}) using  $\overline{\mathbf{A}}$ and $\overline{\mathbf{y}}$, which are normalized $\bA$ and $\by$, respectively.
			\ENDFOR
\STATE Compute $s_k$ in \eqref{eq:cand_score} for $k=1,\dots,K_{\max}-L$, and  
select the optimal candidate with sparsity $k^*$ in \eqref{eq_thresh}. 

   \STATE Obtain $\bc^*(k^*)$ by least square regression using partial columns of $\bA$
   \begin{align*}
	\min_{\mathbf{c}}&\|\mathbf{A}\mathbf{c}-\mathbf{y}\|^2_2~\mbox{ subject to } \text{supp}(\mathbf{c}) =\text{supp}(\overline{\bc}^*(k^*)), 
\end{align*}
or simply rescale $\overline{\bc}^*(k^*)$ according to the norms of columns of $\bA$ and $\by$.
			\RETURN A PDE model specified by $\bc^*(k^*)$.
		\end{algorithmic}
		\caption{The Proposed GP-IDENT  Algorithm}\label{alg_Prop}
	\end{algorithm}

\subsection{Setup: construction of the feature system}\label{ss:feature}  In the first step, we set up a feature system for feature terms identification and coefficients reconstruction.
To simplify the notations, we focus on one-dimensional spacial domain in the  description.

Consider an evolution PDE in~\eqref{eq_PDE_general}  on the spatial and temporal domain $\mathbb{S}^1\times[0,T]$ with a periodic boundary condition in space. Denote a collection of noisy observations of its solution trajectory by $\mathcal{D}=\{U(x_i,t_n)=u(x_i,t_n)+\varepsilon_{i,n}$, $i=1,\dots, I$, $n = 1,\dots,N\}$. Here $\varepsilon_{i,n}$ is the data noise. We assume that  the underlying PDE is in the form of~\eqref{eq_PDE_spec1}, i.e., it is a linear combination of features, e.g., $u_x$ and $uu_x$, contained in an over-complete dictionary $\mathcal{G}=\{f_g:\mathbb{S}^1\times[0,T_{\max}]\mapsto\mathbb{R}\}_{g=1}^G$ with coefficients that may depend on space and time. Note that~\eqref{eq_PDE_spec1} can represent nonlinear PDEs if $\mathcal{G}$ includes nonlinear features.
Let $\{B_m(x,t)\}_{m=1}^M$ be a set of bases, and denote $\mathcal{H}_M = {\rm Span}(\{B_m(x,t)\}_{m=1}^M)$ as a hypothesis space.  We first approximate each variable coefficient $C_g(x_i,t_n)$ by an expansion of the basis elements such that $$C_g(x_i,t_n)\approx \sum_{m=1}^Mc_{g,m}B_m(x_i,t_n) \in \mathcal{H}_M$$ with constant coefficients $c_{g,m}\in \mathbb{R}$  for $g=1,2,\dots,G$.
Then each term in (\ref{eq_PDE_spec1}) is represented as
		\begin{align}
C_g(x_i,t_n)f_g(x_i,t_n)\approx \sum_{m=1}^Mc_{g,m} B_m(x_i,t_n)f_g(x_i,t_n)\;,i=1,\dots,I,~ n= 1,\dots, N. \label{eq_Ck.1}
		\end{align}
		Since the exact value of $f_g(x_i,t_n)$ is  unknown, we approximate it by the empirical counterpart  $\widehat{f}_g(x_i,t_n)$ estimated from the given data $\mathcal{D}$, which is detailed in Section~\ref{Sec:implementation}. We express \eqref{eq_Ck.1} in the matrix form:
		\begin{align}
		  C_g(x_i,t_n)f_g(x_i,t_n)\approx \bff^{\top}_g(i,n)\bc_g,\label{eq_approx}
		\end{align}
		where 
		\begin{align}
		    \bff^{\top}_g(i,n)=\begin{bmatrix}
		        \widehat{f}_g(x_i,t_n)B_1(x_i,t_n) &\cdots &\widehat{f}_g(x_i,t_n)B_M(x_i,t_n)
		    \end{bmatrix} \in \mathbb{R}^{M}
		\end{align}
		and $\bc_g=\begin{bmatrix}
		    c_{g,1}& c_{g,2} & ... &c_{g,M}
		\end{bmatrix}^{\top}\in\mathbb{R}^M$.
Define the $g$-th group feature as
		$$\bF_g=\begin{bmatrix}
		    \bff_g(1,1)&\bff_g(2,1) & \cdots & \bff_g(I,N)
		\end{bmatrix}^{\top}\in \mathbb{R}^{IN\times M}.
		$$
We concatenate $\{\bF_g\}_{g=1}^G$ to construct the {\bf feature matrix}:
		\begin{align}
		    \bA=\begin{bmatrix}
		        \bF_1 & \bF_2 & \cdots & \bF_G
		    \end{bmatrix}\in \mathbb{R}^{IN\times GM},
		\end{align}
as illustrated in Figure~\ref{fig_workflow} (a). Similarly, we construct $\bc$ from $\bc_g$ via
\begin{align}
		    \bc=\begin{bmatrix}
		        \bc_1^\top & \bc_2^\top & \cdots & \bc_G^\top
		    \end{bmatrix}^\top\in \mathbb{R}^{GM}.
		\end{align}
We approximate $u_t(i,n)$  by its empirical counterpart $\widehat{u}_t(i,n)$ based on the given data $\mathcal{D}$. We define the {\bf feature response} as
		\begin{align}
		\by = 
		\begin{bmatrix}
			\widehat{u}_t(x_1,t_1) &
			\widehat{u}_t(x_2,t_1) &
			\cdots &
			\widehat{u}_t(x_I,t_N)
		\end{bmatrix}^{\top}\in \mathbb{R}^{IN},\label{eq_feature_response}
	\end{align}
and refer to the pair $(\mathbf{A},\by)$ as a \textbf{feature system} derived from the given data $\mathcal{D}$ using dictionary $\mathcal{G}$ and the hypothesis space $\mathcal{H}_M$.

 \subsection{Candidate generation using GPSP}\label{ss:candidate}
In the second step, we generate a sequence of candidate models with distinct levels of sparsity. 
Let $K_{\max}$ be a fixed  integer such that $1\leq K_{\max}\leq G$. For  $k=1,2,\dots, K_{\max}$, we consider
\begin{align}
\min_{\overline{\mathbf{c}}\in\mathbb{R}^{GM}}&\|\overline{\mathbf{A}}\overline{\mathbf{c}}-\overline{\mathbf{y}}\|^2_2~\mbox{ subject to }\|\overline{\mathbf{c}}\|_{\ell_{0,1}} = k, \label{eq_cand}
\end{align}
where $\overline{\bA}$ and $\overline{\by}$ are obtained from $\bA$ and $\by$ by normalizing each column, i.e., the column norms are $1$, and 
\begin{align}\|\mathbf{c}\|_{\ell_{0,1}} :=\left\|\begin{bmatrix}\|\mathbf{c}_1\|_1 & \dots &\|\mathbf{c}_G\|_1\end{bmatrix}\right\|_0\end{align}
represents the number of groups with non-zero coefficients. The constraint enforces group sparsity by explicitly specifying that only $k$ groups of features have nonzero coefficients. The solution of~\eqref{eq_cand} corresponds to a PDE model with exactly $k$ features that best fits the given data. However, due to the $\ell_0$-norm constraint, exactly solving the non-convex and non-differentiable problem~\eqref{eq_cand} is NP-hard~\cite{dai2009subspace}.

We propose 
Group Projected Subspace Pursuit (GPSP) to find a group $k$-sparse vector $\overline{\bc}^*(k)$ for $k=1,2,\dots,K_{\max}$. 
Given a fixed $k$, the proposed GPSP iteratively searches for $k$ groups highly correlated to the residuals in a greedy manner (See Section~\ref{sec_GPSP}). 
For each sparsity level $k$, we denote the index set corresponding to the active group features  by $T(k)\subseteq \{1,\dots, G\}$.  We obtain $K_{\max}$ candidate PDEs whose active features are indexed by  $T(k)$ for each sparsity level  $k=1,\dots, K_{\max}$ respectively.

\subsection{Model selection by  Reduction in Residual (RR)}
\label{ss:rss}

The third step is  to select the optimal model from candidates specified by each sparsity level $k$. We design a new score using the residual sum of squares,
$$
R_k=\|\overline{\bA}\overline{\bc}^*(k)-\overline{\by}\|_2^2,
$$
and compare the reduction of this residual for each $k$ sparsity level.  
Let $L\geq 1$ be a fixed integer. For $k=1,\dots, K_{\max}-L$, we compute the Reduction in Residual (RR) as  
 \begin{align}
 s_k = \frac{R_{k}-R_{k+L}}{LR_1}\;,\;\; k=1,\dots, K_{\max}-L.  \label{eq:cand_score}
 \end{align}
This measures the average reduction of residual error as the sparsity level $k$ increases. A small value in $s_k$ means there is a marginal gain in accuracy as sparsity level gets bigger than $k$. 
Here, using $L=1$ is not reliable:  using GPSP, for each sparsity level $k$ the computation of \eqref{eq_cand} is totally independent, the index set $T(k)$ of the active features for the $k$-th candidate, may not be a subset of $T(k+1)$, i.e., $R_k-R_{k+1}$ may be negative.  By using the average of $L$ in~\eqref{eq:cand_score}, we suppress the impact of fluctuation and improve the stability of model selection. 

When the value $s_k$ is already small, we choose the smallest sparsity $k$, rather than choosing $k$ with the smallest $s_k$.  We introduce a  threshold $\rho>0$, and pick the optimal sparsity as follows:
  \begin{align}
 k^*=\min\{k:1\leq k\leq K_{\max}-L, s_k<\rho\}. \label{eq_thresh}
 \end{align}
This is the smallest sparsity index $k$ where $s_k$ is below $\rho$. The motivation of this criterion is to find the simplest model, where RR does not reduce further by considering more complex models.  For the least square fitting, as more terms are added, the error always reduces, RR helps to keep simplest model being independent to increasing level of complexity with increasing sparsity level $k$.  
We find that GP-IDENT is not sensitive to the choice of $L$ and $\rho$, and we fix  $L=5$ and $\rho =0.015$ in this paper.  
We illustrate the effect of RR with an example in Appendix \ref{Asec_RR}.

\subsection{Reconstruction of the coefficients}\label{ss:coeff}

In the fourth step, we reconstruct the coefficients of the identified PDE. After obtaining the optimal level of sparsity $k^*$ in Step 3, we reconstruct the coefficients $\bc^*(k^*)$ by solving
 \begin{align}
	\min_{\mathbf{c} \in \mathbb{R}^{GM}}&\|\mathbf{A}\mathbf{c}-\mathbf{y}\|^2_2~\mbox{ subject to } \text{supp}\,\mathbf{c} =\text{supp}\,\overline{\bc}^*(k^*), \label{eq_rec}
\end{align}
where we recall that $\overline{\bc}^*(k^*)$ is the approximate solution of~\eqref{eq_cand} given by GPSP with the optimal group sparsity $k^*$ selected in subsection~\ref{ss:rss}. It is equivalent to a least square regression using the group features indexed by $T(k^*)$. Alternatively, we can reconstruct $\bc^*(k^*)$ by properly rescaling $\overline{\bc}^*(k^*)$ by the norms of columns of $\bA$ and $\by$. In particular,  the $m$-th entry  $\bc^*(k^*)$ is equal to the $m$-th entry of $\overline{\bc}^*(k^*)$ divided by the norm of the $m$-th column of $\bA$ then multiplied by the norm of $\by$.

 \section{Group Projected Subspace Pursuit (GPSP)}\label{sec_GPSP}

We propose the Group Projected Subspace Pursuit (GPSP) to generate candidates with $k$ features.
For Group-LASSO (GLASSO)~\cite{yuan2006model} and the grouped version of STRidge, Sequential Grouped Threshold Ridge Regression (SGTR)~\cite{Rudy2019DataDrivenIO}, 
the sparsity level is implicitly controlled by a regularization parameter. 
GPSP allows one to explicitly specify the sparsity level, which makes the generation of the candidate models more efficient. Compared to Block Subspace Pursuit (BSP)~\cite{kamali2013block}, GPSP is numerically more stable when co-linearity occurs, and we show this in 
numerical experiments.

For the simplicity of notation, in this section, we use $\bA$ and $\by$ instead of $\overline{\bA}$ and $\overline{\by}$, as the proposed GPSP is applicable in both cases and the normalization is used to make the algorithm numerically robust. 
  
\subsection{GPSP Algorithm} 
For a fixed level of group sparsity $k\geq 1$,  suppose the set of group indices selected by the $l-1$-th iteration is $T^{l-1}$, and denote as
\begin{align}
	      \mathbf{\by}^{l-1}_r =\resid(\by, \bA_{T^{l-1}})= \mathbf{y}-\text{proj}(\by,\bA_{T^{l-1}})=\by-\bA_{T^{l-1}}\bA_{T^{l-1}}^\dagger\mathbf{y}\label{eq_y_residual}
	  \end{align} 
the residual of fitting the data  using groups specified by indices in $T^{l-1}$. Here $\bA_{T^{l-1}}$ is obtained by concatenating the group features $\{\bF_g\}_{g\in T^{l-1}}$ horizontally. The proposed scheme consists of two stages in each iteration: expanding and shrinking.

\noindent\textbf{[Stage 1] Expand $T^{l-1}$ to $\wT^{l}$.} For the $l$-th iteration, we first compute
 \begin{align}
P(\by_r^{l-1},\bF_g)=\frac{\left|\text{proj}(\by_r^{l-1},\bF_g)^{\top}\by_r^{l-1}\right|}{\|\text{proj}(\by_r^{l-1},\bF_g)\|_2\|\by_r^{l-1}\|_2}\label{eq:PyrPg}
 \end{align}
for $g=1,2,\dots, G$.  Note that $P(\by_r^{l-1},\bF_g)$ measures the correlation between $\by_r^{l-1}$ and its projection to the column space of $\bF_g$.  
We take the union of $T^{l-1}$ with the set of $k$ groups with the highest $k$ values in~\eqref{eq:PyrPg}, and denote the union set as $\wT^{l}$. 

\noindent\textbf{[Stage 2] Shrink $\wT^{l}$ to $T^l$.}
 Let $\mathbf{x}^l_p = \bA_{\wT^l}^\dagger\by$.
We project $\by$ to the column space of $\bA_{\wT^l}$ with decomposition $$\by_p = \text{proj}(\by,\bA_{\wT^l})= \sum_{g\in \wT^l }\bF_g \mathbf{x}^l_p[g],$$ where $\bx^l_p[g]$ is the subvector of $\bx^l_p$ corresponding to the $g$-th group.
   For $g\in\wT^l$, its norm
   $\|\bF_g \mathbf{x}^l_p[g]\|_2$ provides a measure of the importance of the $g$-th group.
   Hence, from $\wT^l$, we keep indices of $k$ most important groups and remove the others. The refined set of indices is $T^{l}$. 
   
After the $l$-th iteration, we compute $\by_{r}^{l} =\resid(\by, \bA_{T^{l}})$. If $\|\by_{r}^{l}\|_2 > \|\by_r^{l-1}\|_2$, we take groups indexed by $T^{l}$ as our final selection; otherwise, we  repeat the procedure described above. We summarize  GPSP scheme in Algorithm~\ref{alg_GPSP}.

 \begin{algorithm}[t]
		\begin{algorithmic}[1]
			\REQUIRE{Feature system $(\mathbf{A}$, $\by)$, specified level of group sparsity $k\geq 1$, maximal number of iterations $\text{Iter}_{\max}\geq 1$.
			}
			\STATE Set $l=0$.
			\STATE Set 
   $T^l = \{k$ largest indices of 
   $P(\by,\bF_g), g=1,2,\dots,G\}$ in \eqref{eq:PyrPg}.
			\STATE Set $\by_r^l = \resid(\by, \bA_{T^l})$ in \eqref{eq_y_residual}, $\bA_{T^l}$ 
   concatenates $\{\bF_g\}_{g\in T^l}$ vertically.	
   
   \FOR{$l=1,\dots,\text{Iter}_{\max}$}
		\STATE$\wT^l = T^{l-1}\cup\{
  k$ largest indices of $P(\by,\bF_g), g=1,2,\dots,G\}$.
\STATE  Compute $\mathbf{x}^l_p = \bA_{\wT^l}^\dagger\by$.
\STATE Set  $T^{l}=\{
   k$ largest indices of 
$\|\bF_g\bx^l_p[g]\|_2,~g\in \wT^l\}$, where $\bx^l_p[g]$ is the subvector of $\bx^l_p$ corresponding to the $g$-th group.
\STATE Compute $\by_r^l = \resid(\by, \bA_{T^l})$.
\IF{$\|\by_r^l\|_2 > \|\by_r^{l-1}\|_2$}
			\STATE Set $T^l = T^{l-1}$ and terminate.
			\ENDIF
			\ENDFOR
			
			\RETURN The optimal group indices $T^l$ and the estimated coefficient $\bA_{T^l}^\dagger\by$
		\end{algorithmic}
		\caption{Group Projected Subspace Pursuit (GPSP) for \eqref{eq_cand}}\label{alg_GPSP}
	\end{algorithm}

\subsection{Related algorithms}
GPSP is closely related to  Subspace Pursuit (SP)~\cite{dai2009subspace} and Block Subspace Pursuit (BSP)~\cite{kamali2013block}.  SP is a greedy algorithm for  sparse regression. 
It iteratively expands the pool of $k$ candidate covariates by considering potential features highly correlated to the residual, then refines the choices by reducing the extended pool back to $k$ covariates by eliminating those with less importance. At each iteration, SP expands $k$ nonzero entries to $2k$ nonzero entries by adding the $k$ indices whose columns are highly correlated with the residual, and then refines the choice by eliminating the $k$ indices with smaller coefficient values. From this perspective, both BSP and GPSP can be regarded as generalizations of SP where the covariates, i.e., individual columns of the system matrix, are replaced by groups of features. However, BSP and GPSP have different interpretations about the correlation between the residual and a feature group.  

\begin{figure}
		\centering
		\begin{tabular}{cc}
			(a)&(b)\\
			\includegraphics[width = 0.4\textwidth]{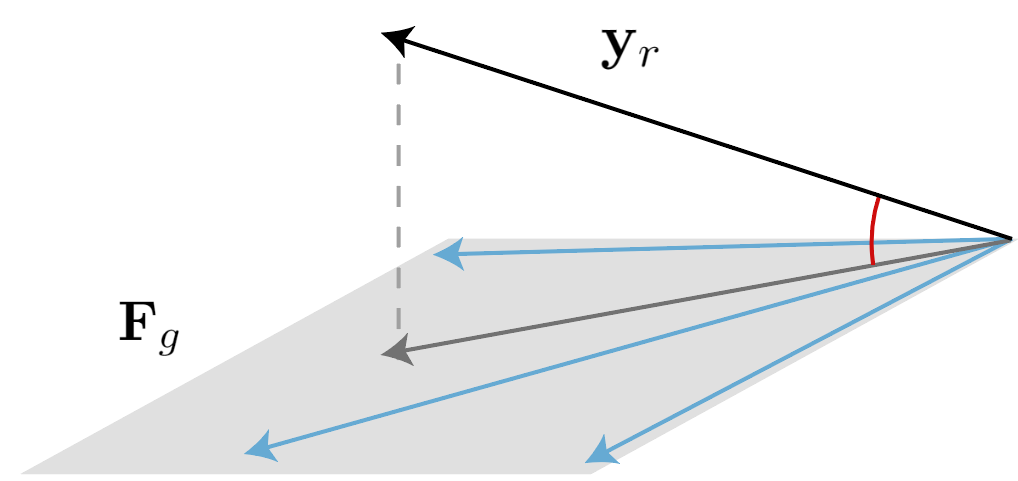}&
			\includegraphics[width=0.4\textwidth]{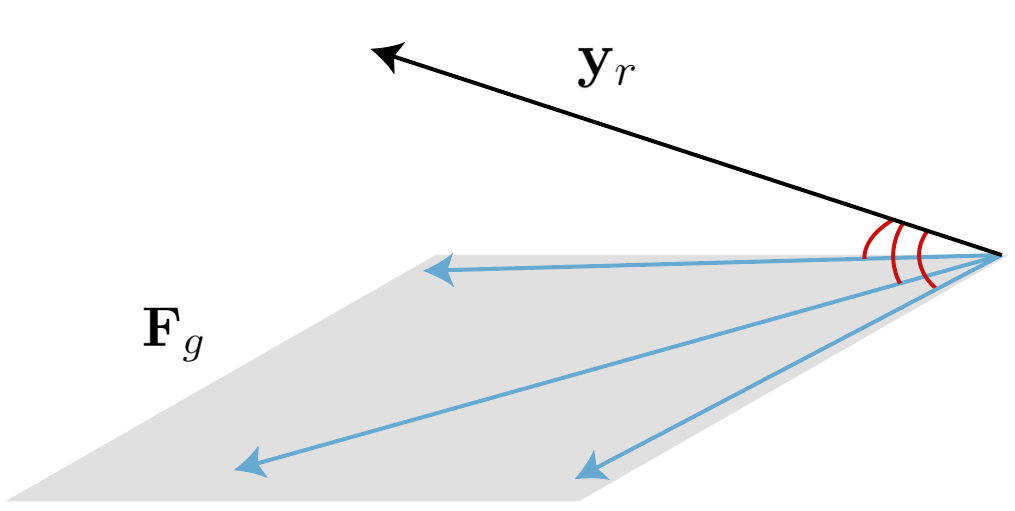}
		\end{tabular}
\caption{An illustrative comparison between GPSP and BSP~\cite{kamali2013block}.  (a) In GPSP, the group's importance is evaluated by the correlation between $\by_r$ and its projection to the column space of $\bF_g$.  (b) In BSP, the group's importance is evaluated by the correlation between the residual  $\by_r$ and the columns (blue arrows) in $\bF_g$. }\label{fig_bsp_gpsp}
\end{figure}

When expanding the pool of candidates from $k$ to $2k$, BSP measures the correlation between the residual $\by_r$ and the $g$-th feature group $\bF_g$ by the $L_2$-norm of the  inner product between  $\by_r$ and  the columns of $\bF_g$, 
\begin{align}
	\|\bF_g^T\by_r\|_2=\sqrt{\sum_{m=1}^M(\bF^\top_g[m]\by_r)^2}\label{eq_bsp_in},
\end{align}
where $\bF_g[m]$ denotes the $m$-th column of the $g$-th feature group. In GPSP, we use the inner product between  $\by_r$ and its projection to the column space of $\bF_g$ to quantify the correlation 
\begin{align}
	P(\by_r,\bF_g) = \frac{\left|(\bF_g\bF_g^\dagger\by_r)^T\by_r\right|}{\|\bF_g\bF_g^\dagger\by_r\|_2\|\by_r\|_2}. \label{eq_gpsp_in}
\end{align}
Comparing~\eqref{eq_bsp_in} with~\eqref{eq_gpsp_in}, we note that GPSP is less sensitive to co-linearity than BSP.  If some columns of $\bF_g$ are co-linear, BSP~\eqref{eq_bsp_in} considers that they all contribute to the correlation between $\by_r$ and $\bF_g$, whereas GPSP~\eqref{eq_gpsp_in} ignores the co-linear columns as they are redundant when representing the information contained in the group. See Figure~\ref{fig_bsp_gpsp} for an illustration. Notice that if $\bF_b$ only has one column, both~\eqref{eq_bsp_in} and~\eqref{eq_gpsp_in} are identical to SP. 

When reducing the expanded pool of candidates of size $2k$ to $k$, BSP  keeps the $k$ groups whose reconstructed coefficients have the largest magnitudes, whereas  GPSP  uses each group's contribution measured by the norm of the response vector.
Table~\ref{tab_BSP_GPSP} summarizes the differences between BSP and GPSP. 
\begin{table}
	\centering
	\begin{tabular}{c|c|c}
		\hline
		\textbf{Criterion}&\textbf{BSP}~\cite{kamali2013block}&\textbf{GPSP} (Proposed)\\\hline
		\Gape[0.1cm][0.1cm]{Expand} &$	\|\bF_g^\top\mathbf{y}_r^{l-1}\|_2$&$\left|\bF_g\bF_g^\dagger\by_r^{l-1})^\top\by_r^{l-1}\right|/(\|\bF_g\bF_g^\dagger\by_r^{l-1}\|_2\|\by_r^{l-1}\|)$\\\hline
		\Gape[0.2cm][0.2cm]{Shrink }&$\mathbf{A}_{\wT^l}^\dagger\mathbf{y}$&$\bA_{\wT^l}\bA_{\wT^l}^\dagger\mathbf{y}$\\\hline
	\end{tabular}
\caption{Comparison {of two stages} in BSP and GPSP. In [Stage 1] Expand (the first row), BSP chooses the groups $\bF_g$ whose columns are highly correlated with the residual $\by_r^{l-1}$, whereas GPSP chooses the groups whose column spaces are close to the residual. {In [Stage 2] Shrink (the second row),}  BSP selects the groups with large coefficients, while GPSP selects the groups whose projected residual is significant. } \label{tab_BSP_GPSP}
\end{table}
In general, GPSP is better suited for identifying PDEs with varying coefficients which are approximated by a basis expansion. As we allow the coefficients to vary both in space and time, some columns in the feature matrix can be highly correlated. We observe that GPSP is more effective than BSP when some columns within the same group are highly correlated.
We illustrate the effect with the transport equation with constant speed $a\neq 0$, 
$ u_t(x,t) = au_x(x,t)\label{eq:amb_eq1}$
in Appendix \ref{sec_GPSP_BSP_PDE}.
We also numerically justify these in Section~\ref{sec_exp}. 

\section{Numerical Implementation Details} \label{Sec:implementation}
In this section, we present computational details for B-spine set-up and details of SDD used in this paper. 

\subsection{Approximation of varying coefficients by B-splines}\label{ss_basis}
For some fixed integer $M\geq 1$, we define $\mathcal{H}_M=\{\sum_{m=1}^Mc_mB_m:c_m\in\mathbb{R}\}$ as our hypothesis space, where the basis function $B_m\in\mathbb{S}^1\times[0,T_{\max}]\mapsto \mathbb{R}$ is compactly supported and
$\sum_{m=1}^M B_m(x,t)=1$ for all $(x,t) \in \mathbb{S}^1\times[0,T_{\max}]$. The function space $\mathcal{H}_M$ is used to approximate the varying coefficients in the PDE, and we use the basis functions $B_m$'s given by B-splines~\cite{piegl1996nurbs}.

Without loss of generality, we consider $[0,1]$ as the spacial domain of interest. For a fixed integer $p\geq 1$, we consider a uniform knot sequence $0= z_0<z_1<\cdots<z_l= 1$ for some $l\geq p$. Denote the knot spacing by $\Delta z$. 
The $n$-th B-spline basis function $b^p_n$ of order $p$ is constructed according to the Cox-de Boor recursion formula~\cite{de1986b}
\begin{align}
b^0_n(z) &= \begin{cases}
z&\text{if}~ z_n\leq z<z_{n+1},\\
0&\text{otherwise},
\end{cases}\\
	b^p_n(z) &=\left(\frac{z-z_n}{z_{n+q}-z_n}\right)b_n^{p-1}(z) + \left(\frac{z_{n+p+1}-z}{z_{n+p+1}-z_{n+1}}\right) b_{n+1}^{p-1}(z).
\end{align}
for $0\leq n\leq l-p-1$.
We note that $b_n^p$ is non-zero on $[z_n, z_{n+p+1})$, and there are at most $p+1$ non-zero basis functions over any interval $[z_n,z_{n+1})$. Suppose the knot spacing is $\Delta z$, depending on different boundary conditions for the functions to be approximated, we supplement $\{b_n^p\}_{n=0}^{l-p-1}$ with more basis functions. For this purpose, it is convenient to uniformly extend the knot sequence to infinity $\cdots<z_{-2}<z_{-1}<z_0<\cdots<z_l<z_{l+1}<\cdots$ where $b_n^p$ is defined  for $n\in\mathbb{Z}$. 
\begin{itemize}
\item \textbf{Periodic boundary condition.} Add $p$ functions  $\widetilde{b}_{n}^p$ for $n=-p,-p+1,\dots,-1$ defined as
\begin{align}
\widetilde{b}_{n}^p(z) = \begin{cases}
b_n^p(z)&\text{if}~0\leq z< (n+p+1)\Delta z,\\
b_n^p(z-1)&\text{if}~1+n\Delta z\leq z\leq  1,\\
0&\text{otherwise}.
\end{cases}
\end{align}
\item \textbf{Neumann boundary condition.} Add two  functions 
\begin{align}
&b_L^p(z) =  \begin{cases}
\sum_{n=-p}^{-1} b_n^p(z)&\text{if}~0\leq z< p\Delta z,\\
0&\text{otherwise},
\end{cases}\\
&b_R^p(z) =  \begin{cases}
\sum_{n=l-p}^{l-1} b_n^p(z)&\text{if}~1-p\Delta z \leq z \leq 1,\\
0&\text{otherwise}.
\end{cases}
\end{align}
\end{itemize}
It is easy to check that $\sum_{n=0}^{l-p-1}b_n^{p}(z)+\sum_{n=-p}^{-1}\widetilde{b}_n^p (z)= 1$ and  $\sum_{n=0}^{m-p-1}b_n^{p-1}(z)+b_L^p(z)+b_R^p(z) =1$ when  $z\in [0,1]$, and $\{b_n^{p}\}_{n=0}^{l-p-1}\cup\{\widetilde{b}_n^p\}_{n=-p}^{-1}$ serve as a set of B-spline basis functions of order $p$ for $\mathbb{S}^1$.

In this paper, we assume periodic boundary condition in the space. In the time direction,
we assume the Neumann boundary condition that the underlying coefficients do not have significant changes at the first nor the last moment of the observation. Suppose $\{b_{m_1}(x)\}_{m_1=1}^{M_1}$ is a set of B-spline bases constructed for $\mathbb{S}^1$, and another set $\{b_{m_2}(t)\}_{m_2=1}^{M_2}$ is constructed for $[0,T_{\max}]$ with supplementary elements for the  Neumann boundary condition. We obtain a set of B-spline bases on the spatio-temporal domain $\mathbb{S}^1\times [0,T_{\max}]$ by taking tensor products, that is, 
$$B_m(x,t)\in \{b_{m_1}(x)b_{m_2}(t):m_1=1,\dots,M_1,~m_2=1,\dots,M_2\}$$
for $m=1,2,\dots,M$, where $M=M_1M_2$.

 \subsection{SDD for robust feature approximation}\label{ss:sdd}
To robustly approximate $f_g$ (respectively $u_t$)  with $\widehat{f}_g$ (respectively $\widehat{u}_t$) using noisy observations of $u$~\eqref{eq_approx}, we suppress the noise amplification during the process of numerical differentiation.  We apply the Successively Denoised Differentiation (SDD)~\cite{he2022robust}, which approximates
$\partial_x^n\partial_t^m u(i,j)$ for any integers $m,n\geq 0$ by
$$(\mathcal{S}_xD_x)^n(\mathcal{S}_tD_t)^m\mathcal{S}_x\mathcal{S}_tU(i,j)$$
where $\mathcal{S}_x$ and $\mathcal{S}_t$ are 1-D smoothing operators along space and time respectively, $D_x$ and $D_t$ are numerical differentiation operators with respect to space and time respectively, and $(\cdot)^m$ means applying the operator repeatedly for $m$ times. 

In this paper, we assume that the grid is uniform  with step size $\Delta x>0$ in space and  $\Delta t>0$ in time. We use 5-point-central difference for both $D_x$ and $D_t$
\begin{align*}
D_xU(x_i,t_n) = \frac{-U(x_{i+2},t_n)+8U(x_{i+1},t_n)- 8U(x_{i-1},t_n)+U(x_{i-2},t_n)}{12\Delta x}\;
\end{align*}
and similarly for $D_t$. Here periodic boundary condition is applied for the space, and Neumann boundary condition is applied for the time. To reduce the influence of the approximation errors near boundary, we only use the interior data for feature construction. As for the smoothing operator in time and space, we use the Savitzy-Golay filter~\cite{savitzky1964smoothing}, which is a convolution version of the local polynomial fitting. For example, when the boundary condition is periodic, the spacial smoothing operator with the Savitzy-Golay filter is
\begin{equation} \label{eq_sdd-num}
\mathcal{S}_xU(i,n) = \sum_{l=\frac{1-w}{2}}^{\frac{w-1}{2}}W_lU(i+l,n),
\end{equation}
where the integer $w\geq 1$ is the window size,  the convolution weights $W_l$ are derived by fitting  local data using degree $q$ polynomials for some integer $0\leq q< w$, and they are tabulated in~\cite{savitzky1964smoothing}. This filter is available, e.g., using \texttt{savgol\_filter} from the \texttt{scipy} package in Python. In the following numerical section, we use the notation such as SDD-15 to represent using SDD with a window width $w = 15$ in (\ref{eq_sdd-num}).   In this paper, we find that more accurate coefficient reconstruction is obtained if $\partial_t$ is approximated by $D_t\mathcal{S}_t$ without the second smoothing, thus we modify SDD as such in our experiments.

\section{Numerical Experiments}
\label{sec_exp}
We next present numerical experiments to justify the effectiveness of GP-IDENT and compare it with the state-of-the-art identification methods for varying coefficient PDEs: GLASSO~\cite{yuan2006model}, SGTR~\cite{Rudy2019DataDrivenIO}, and rSGTR~\cite{li2020robust}\footnote{For GLASSO~\cite{yuan2006model} and SGTR~\cite{Rudy2019DataDrivenIO}, we used the code is available at \url{https://github.com/snagcliffs/parametric-discovery}; and for rSGTR, \url{https://github.com/junli2019/Robust-Discovery-of-PDEs}}.  We note that in~\cite{li2020robust}, DLrSR was proposed to handle sparse noise added to the measurements in a linear system.   We also compare GP-IDENT with BSP-IDENT, where GPSP is replaced by BSP~\cite{kamali2013block} in GR-IDENT. To show the effectiveness of GP-IDENT, we test it on various types of  equations~\cite{rudy2017data,Rudy2019DataDrivenIO,li2020robust,he2022robust} listed in Table~\ref{tab_PDE_exp}.  For the PDE examples, we generate the solution data by the spectral method analogous to~\cite{Rudy2019DataDrivenIO}. The equation is discretized in space, where the partial derivatives are computed using Fast Fourier Transform (FFT), then the solution is obtained by integrating in time using LSODA~\cite{LSODA}. As for examples of PDE systems, i.e., Schr\"{o}dinger and Nonlinear Schr\"{o}dinger equations, we generate the data by implicit-explicit finite difference methods where second order partial derivatives are treated implicitly, and the zero-th order terms are treated explicitly.

\begin{table}[t]
	\centering
	\begin{tabular}{l|l}
		\hline
		\textbf{PDE}& \textbf{Model}\\\hline
  Advection diffusion equation& $u_t = \partial_x(a(x)u)+bu_{xx}$\\\hline
  Fisher's equation& $u_t = bu_{xx}+a(t)u(1-u)$\\\hline
  Viscous Burgers equation& $u_t=a(x,t)uu_x+b(t)u_{xx}$\\\hline
  Korteweg–De Vries (KdV) equation&$u_t=a(x,t)uu_x+b(x,t)u_{xx}$\\\hline
  Kuramoto–Sivashinsky (KS) equation& $u_t = a(x)uu_x+b(x,t)u_{xx}+c(x,t)u_{xxxx}$\\\hline
  Schr\"{o}dinger equation& $iu_t=bu_{xx}+a(x,t)u$\\\hline
  Nonlinear Schr\"{o}dinger (NLS) equation& $iu_t=bu_{xx}+a(x,t)|u|^2u$\\\hline
	\end{tabular}
	\caption{A list of PDEs tested in Section~\ref{sec_exp}. }\label{tab_PDE_exp}
\end{table}

We consider data with $p\%$ Gaussian noise.  The noisy data takes the form ${U}(x_i,t_n)$ = $u(x_i,t_n)+\varepsilon_{i,n}$ for $i=1,\dots,I$, $n=1,\dots,N$, with Gaussian noise 
$ \varepsilon_{i,n}\sim \mathcal{N}(0,\sigma^2)$, $i=1,\dots,I$, $n=1,\dots,N
$,  
where
\begin{align}
\sigma = p\%\times\text{std}(\{u(x_i,t_n)~|~i=1,\dots,I,n=1,\dots,N\}).
\end{align}
Here $\text{std}(\cdot)$ stands for the standard deviation of a collection of data.

To evaluate the reconstruction accuracy, we calculate  the discrete relative $L_1$-error to measure the coefficient error:
\begin{align}
e(C_g) = \frac{\sum_{i=1}^I\sum_{n=1}^N\left|\widehat{C}_g(x_i,t_n)-C_g(x_i,t_n)\right|}{\sum_{i=1}^I\sum_{n=1}^N\left|C_g(x_i,t_n)\right|}\times 100\%\label{eq:error}
\end{align}
where $\widehat{C}_g$ is the reconstruction of $C_g$.  To quantify the {the coefficient support } identification accuracy, we use the Jaccard index~\cite{jaccard1912distribution} defined as
\begin{align}
J(\widehat{T},T^*) =\frac{|\widehat{T}\cap T^*|}{|\widehat{T}\cup T^*|},
\end{align}
where $\widehat{T}$ denotes the group index set  in the identified model,  $T^*$ is the group index set in the true equation, and $|\cdot|$ gives the number of elements in the set.  Note that $J(\widehat{T},T^*)=1$ if and only if $\widehat{T}=T^*$, i.e., the underlying model is exactly identified. 

For hyper-parameters, we fix $K_{\max}=15$, $\rho=0.015$, and $L=5$ in all experiments. Our default dictionary contains 56 terms including all partial derivatives of $u$ up to order $4$ and the products of no more than of $3$ features. Our experiments in subsection \ref{sec_adv} are performed on larger dictionaries for comparisons. 

\subsection{GP-IDENT results on PDEs with space and time varying  coefficients}
\begin{table}[t!]
\small
\centering
\begin{tabular}{c|c|c}
\hline
Trajectory & Equation &   Coef. error no noise,   $1\%$ noise   \\\hline
KdV \parbox[c]{4.7em}{
\includegraphics[width=0.6in]{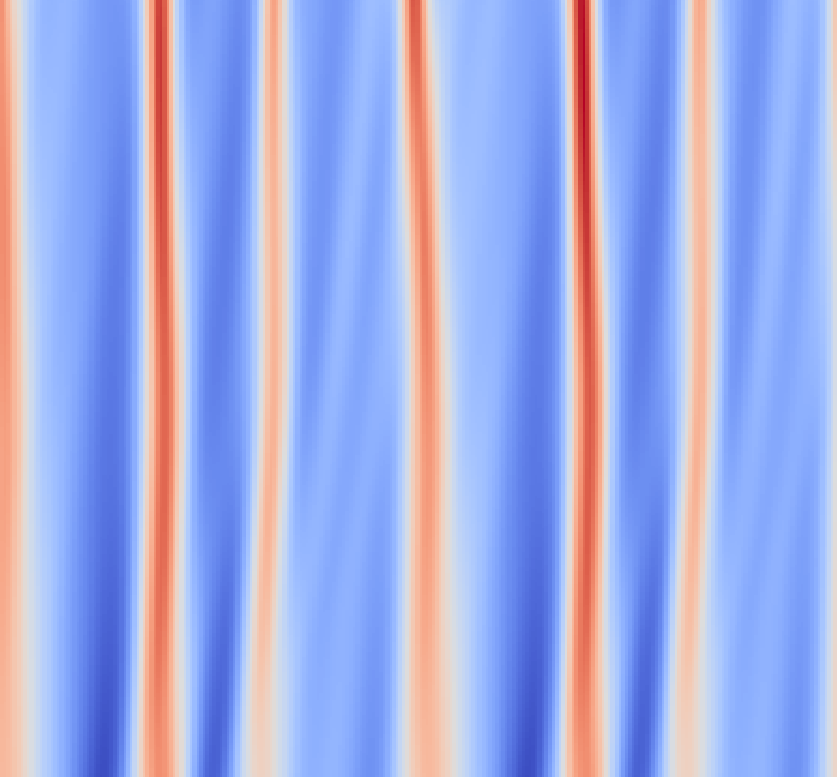}}
&$u_t = a(x,t)\, uu_x + b(x,t)\, u_{xxx}$ & $\begin{aligned}
&u_{xxx}&:4.09\%, \;\; &20.49\pm0.16\%\\
&uu_{x}&:0.54\%, \;\;  &20.37\pm0.18\% 
\end{aligned}$\\\hline
KS~~ \parbox[c]{4.7em}{
      \includegraphics[width=0.6in]{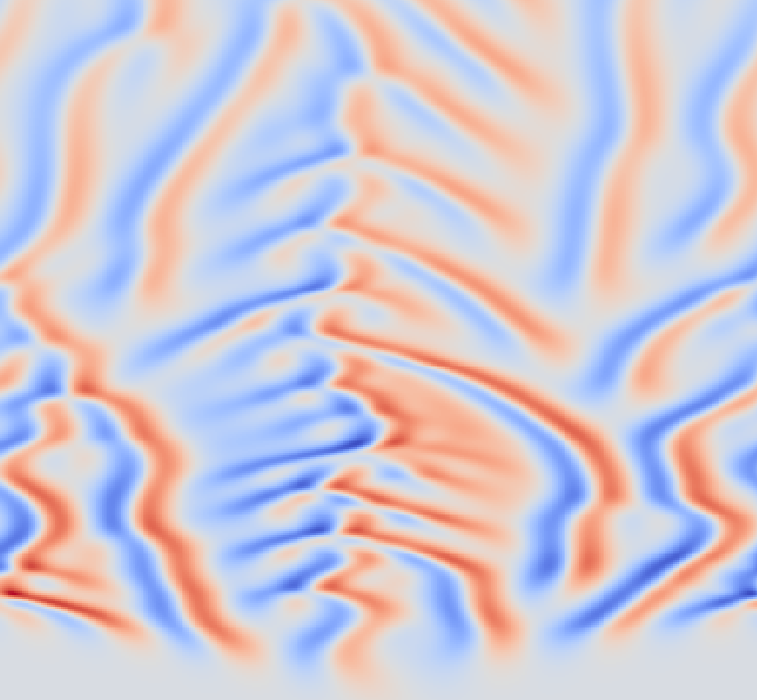}}&
      $\begin{aligned}
u_t = a(x)\,uu_x &+ b(x,t)\,u_{xx}\\
&+c(x,t)\,\partial_x^4u
      \end{aligned}
      $
      &$\begin{aligned}
&u_{xx}&:2.03\%, \;\; 19.21\pm0.31\%^\dagger\\
&\partial_x^4u&:2.12\%,\;\; 18.92\pm0.30\%^\dagger\\
&uu_x&:1.05\%,\;\; 25.61\pm0.21\%^\dagger  
      \end{aligned}$\\\hline
Sch~ \parbox[c]{4.7em}{
      \includegraphics[width=0.6in]{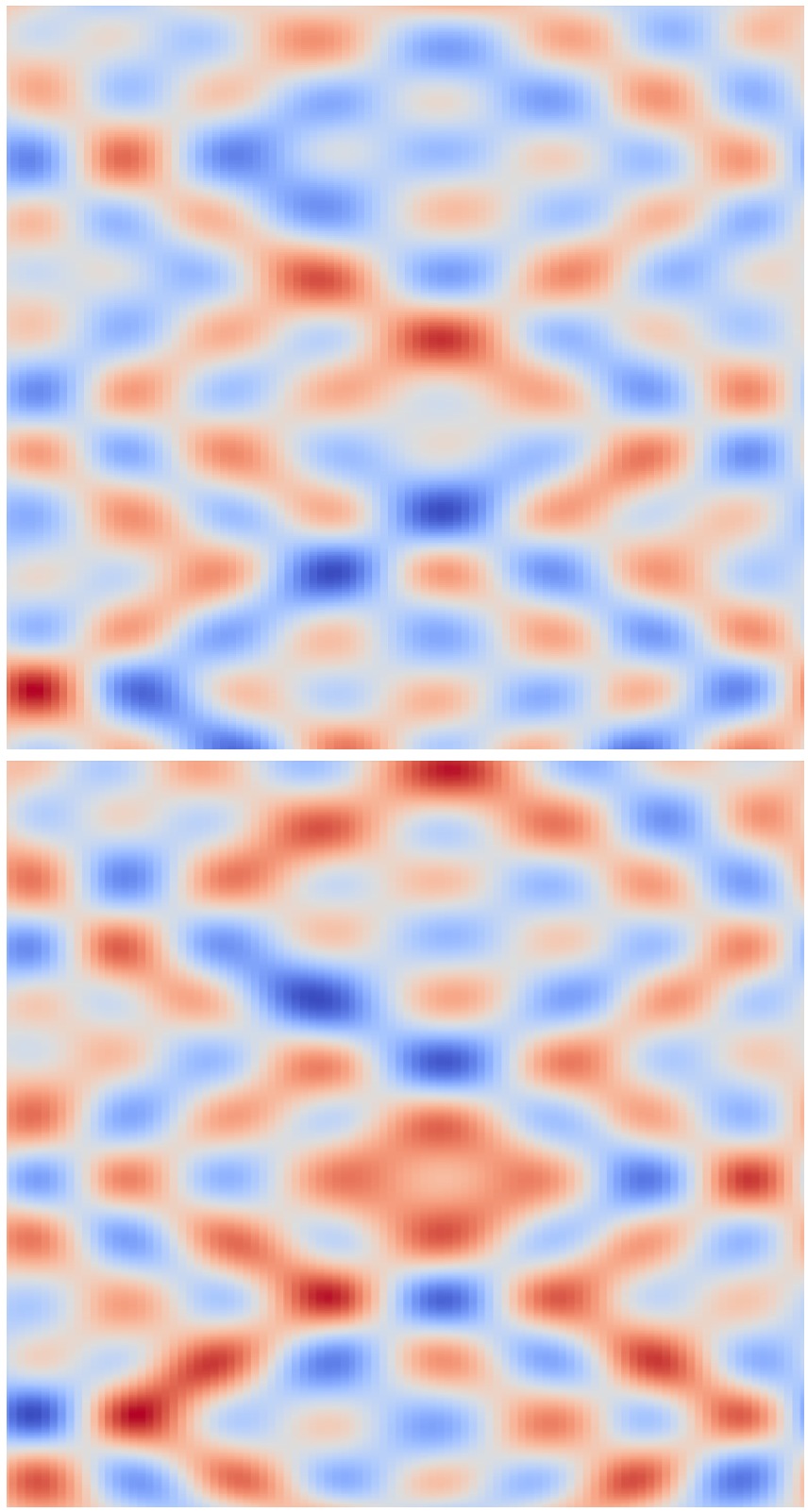}}&$iu_t = 0.5u_{xx}+a(x,t)u$&$\begin{aligned}
&v&:3.80\%,\;\; 5.34\pm0.12\%\\
&w&:3.93\%,\;\; 4.14\pm0.11\%\\
&v_{xx}&:0.78\%, \;\; 0.87\pm0.03\%\\
&w_{xx}&:0.81\%,\;\; 0.73\pm0.03\%\\
&& v=\text{Re}(u), w=\text{Im}(u)
      \end{aligned}$\\\hline
NLS \parbox[c]{4.7em}{
      \includegraphics[width=0.6in]{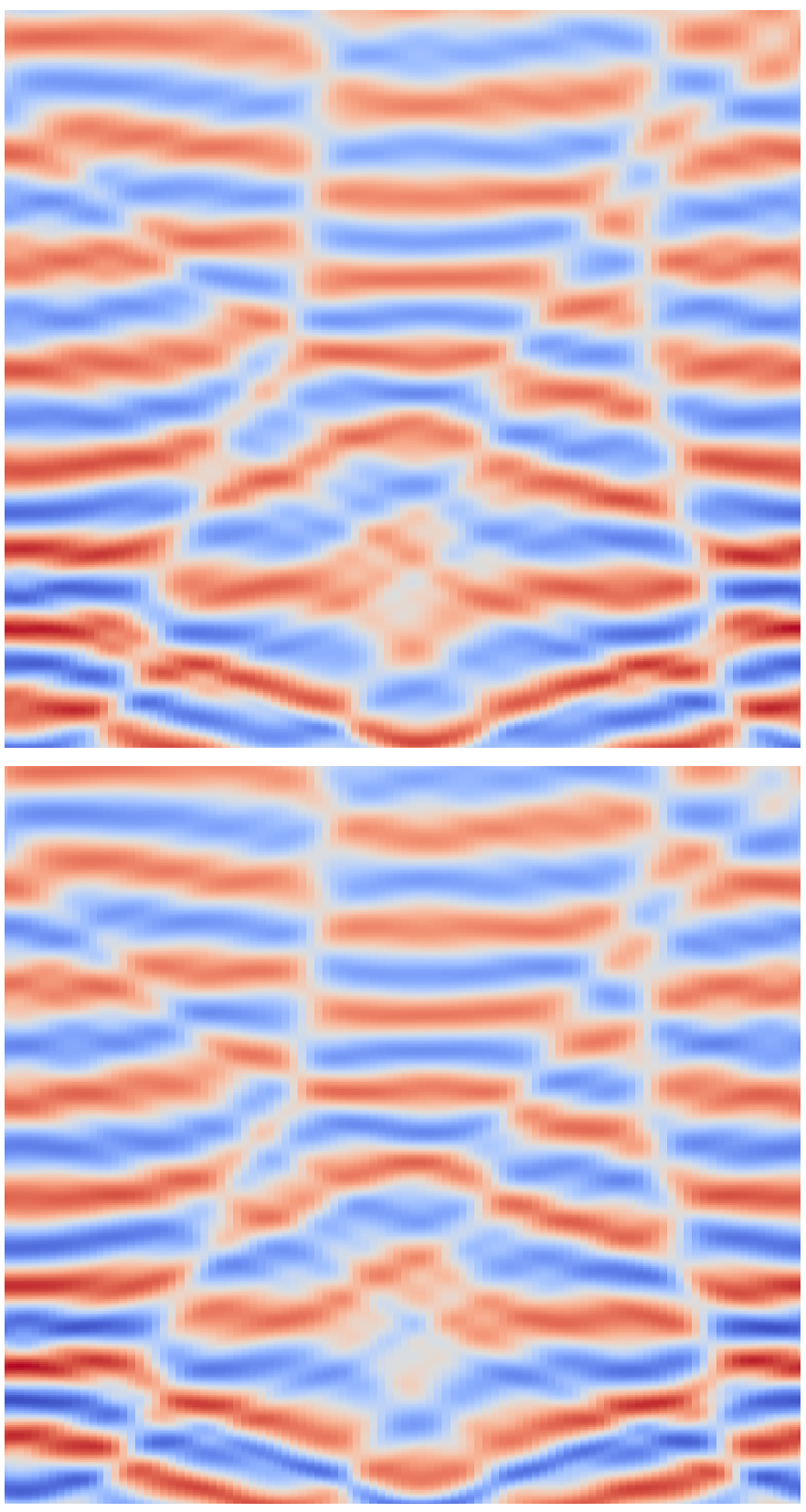}}&$iu_t = -0.5u_{xx}+a(x,t)|u|^2u$&$\begin{aligned}
&v_{xx}&: 1.74\%,\;\; 21.48\pm0.05\%\\
&w_{xx}&: 1.97\%,\;\; 21.16\pm0.04\%\\
&v^3&: 0.39\%,\;\; 2.67\pm0.01\%\\
&v^2w&:0.40\%,\;\; 2.76\pm0.04\%\\
&vw^2&:0.50\%,\;\; 2.61\pm0.02\%\\
&w^3&:0.36\%,\;\; 2.44\pm0.01\%\\
&&v =\text{Re}(u), w=\text{Im}(u)
      \end{aligned}$\\\hline
\end{tabular}
\caption{GP-IDENT results for equations and systems with space and time varying coefficients.
The first column shows the solution trajectory for each equation. For the Schr\"{o}dinger and NLS equations, the real and imaginary components of $u$ are plotted, respectively.  The second column shows the identified equations, whose features match the respective true features. The third column, we report the relative $L_1$ coefficient  errors~\eqref{eq:error} for the identified features. For the noisy case ($1\%$ noise), we conduct 10 independent experiments, and record the mean and standard deviation of the errors. $^\dagger$To identify KS equation from noisy data, we used $\rho=0.05$ for model selection.}\label{tab_more_exp}
\end{table}

We experiment on several PDEs with space and time dependent  coefficients, including the KdV equation, KS equation, Schr\"{o}dinger equation (Sch), and Nonlinear Schr\"{o}dinger (NLS) equation.  We note that Sch and NLS equations can be regarded as PDE systems for the real and imaginary components of a complex system. For the KdV and KS equations, we use the default dictionary containing 56 terms. As for the PDE systems (Sch and NLS), we use the dictionary containing linear features of  partial derivatives of the real and imaginary components up to order $3$, and the products up to $3$ terms,  leading to a total of $165$ features. Table~\ref{tab_more_exp} shows the trajectories, equations, and the coefficient reconstruction errors~\eqref{eq:error}  with clean and noisy data.  We present the details of these experiment settings including the coefficients, grid, number of bases, and window size for SDD in Appendix \ref{Asec_sdd} Table~\ref{tab_more_exp_info}.
For the KS equation, a different threshold $\rho=0.05$ is used.

\subsection{Viscous Burgers' equation with space-time dependent coefficients}

\begin{figure}[t!]
		\centering
		\begin{tabular}{ccc}
			(a)&(b)&(c)\\
			\includegraphics[width = 0.3\textwidth]{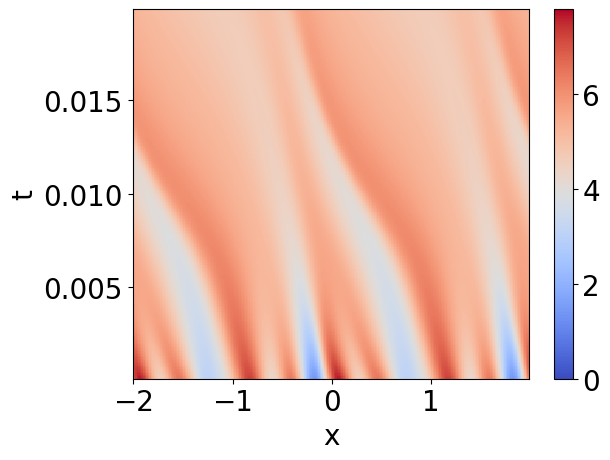}&
			\includegraphics[width=0.3\textwidth]{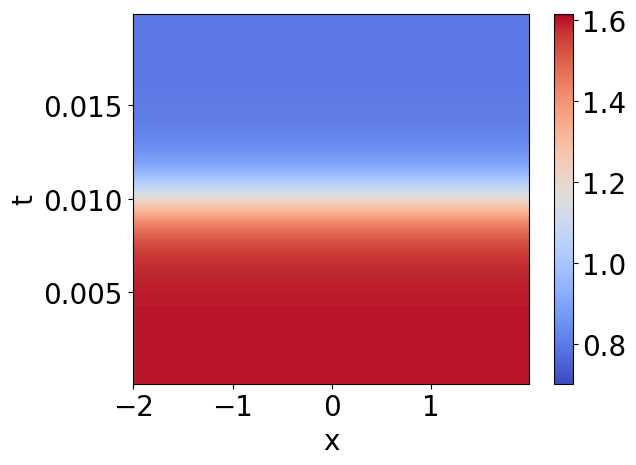}
   &
			\includegraphics[width=0.3\textwidth]{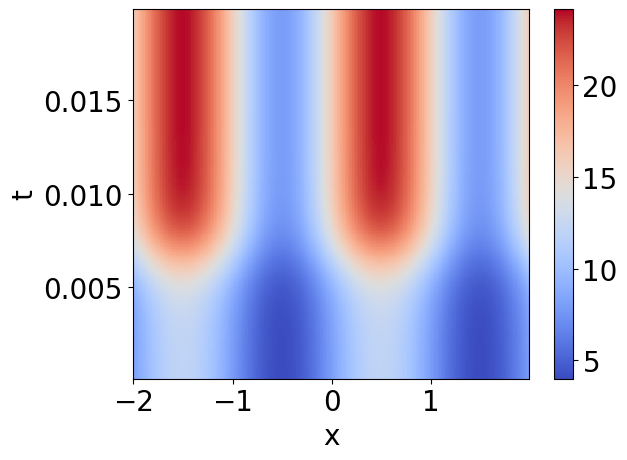}\\
   (d)&(e)&(f)\\
			\includegraphics[width = 0.3\textwidth]{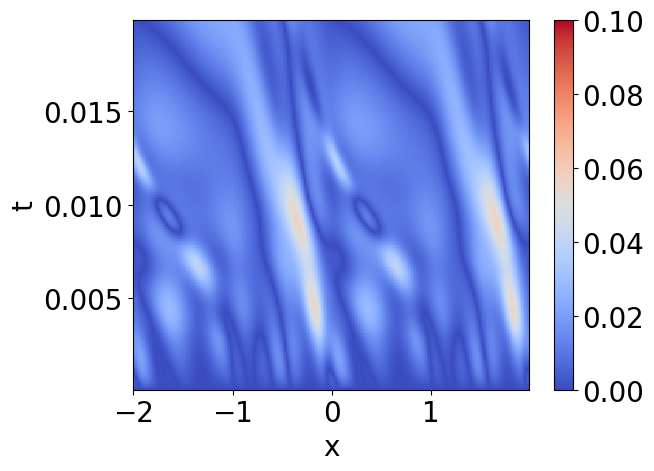}&
			\includegraphics[width=0.3\textwidth]{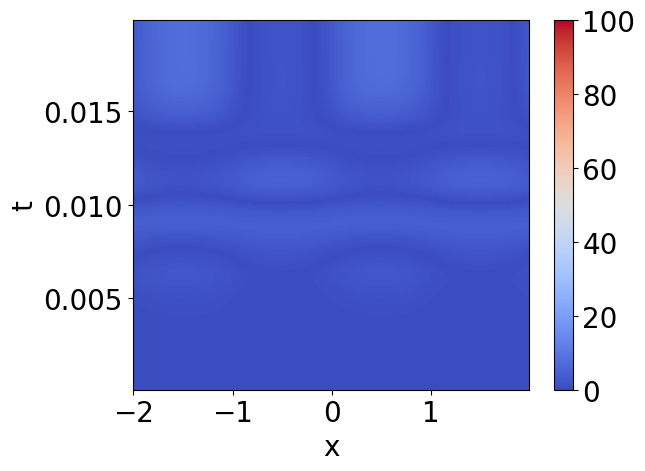}
   &
			\includegraphics[width=0.3\textwidth]{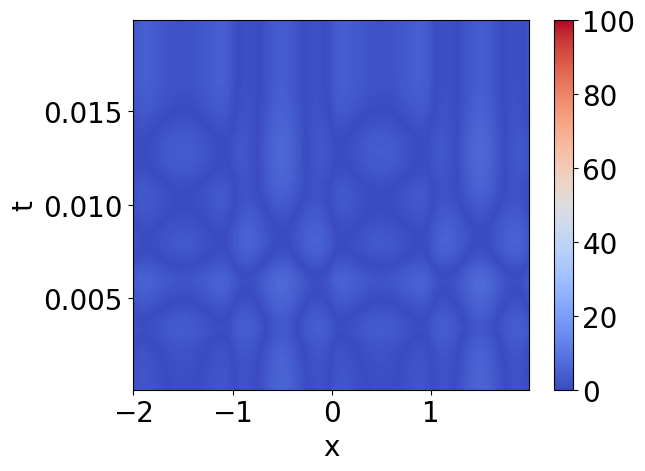}\\
   (g)&(h)&(i)\\
   \includegraphics[width = 0.3\textwidth]{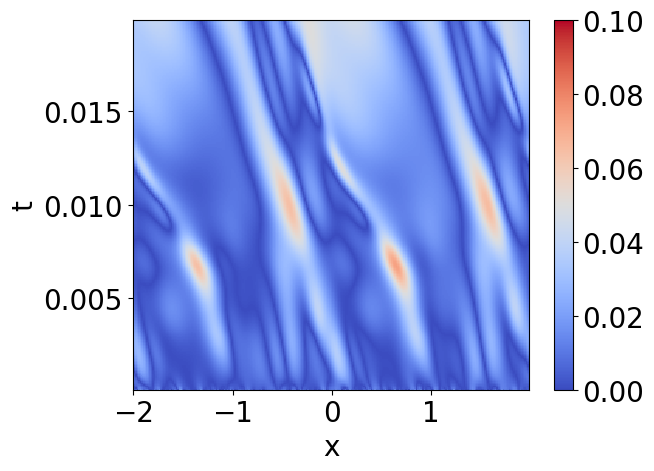}&
			\includegraphics[width=0.3\textwidth]{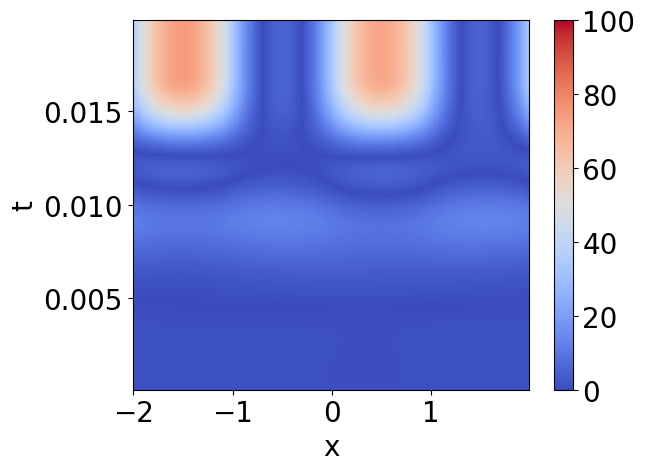}
   &
			\includegraphics[width=0.3\textwidth]{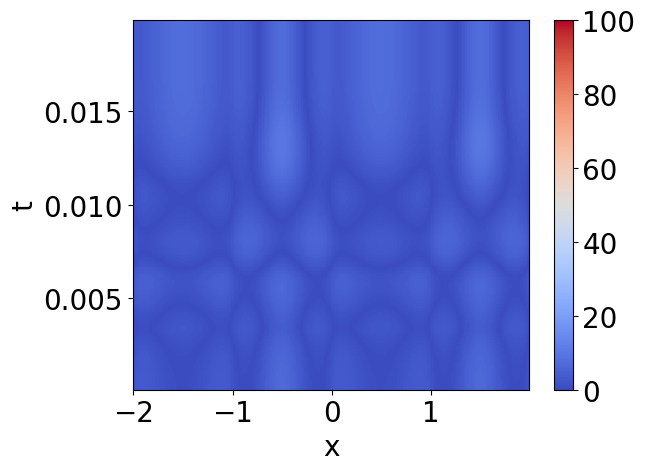}\\
		\end{tabular}
\caption{GP-IDENT result on viscous Burgers equation~\eqref{eq_exp1}. (a) The true clean trajectory,  
(b) the true coefficient for $u_{xx}$, and (c) the true coefficient for $uu_x$.  The second row shows GP-IDENT for clean data: (d) the absolute error of  simulation from the identified model, (e) the relative $L_1$ coefficient error in percentage of the coefficient for $u_{xx}$, and (f) the relative $L_1$ coefficient error in percentage for $uu_x$.  
The third row shows GP-IDENT for the data with 2\% noise: (g) the absolute error of  simulation from the identified model, (h) the relative $L_1$ coefficient error in percentage of the coefficient for $u_{xx}$,  and (i) the relative $L_1$-error in percentage of the coefficient for $uu_x$.}\label{fig_general}
	\end{figure}

Consider the following viscous Burgers' equation
\begin{align}
		u_t(x,t) &= a(x,t)u(x,t) u_x(x,t)+b(t)u_{xx}(x,t)\;,~x\in[-2,2), t\in (0,0.02]\label{eq_exp1}
	\end{align}
with the initial condition 
\begin{align}
		u(x,0)=&\sin(\pi(2x -0.1)) + \cos(\pi(5x -0.2)) + \cos(\pi(3x - 0.3))\cos(\pi(x + 0.1)) \nonumber\\
  &+
         \sin(\pi (4x + 0.5)) + 5\label{eq_exp1_init}
	\end{align}
and space-time dependent  coefficients
\begin{align}
a(x,t)=4\left(1 + \tau_{+}\left(t,10,\frac{0.02}{3}\right)\right) (2+\sin(\pi x))\;,~b(t)=0.8\left(1+\tau_{-}\left(t,10,\frac{0.02}{2}\right)\right).
\end{align}
We numerically solve it on a $256\times256$ grid. Figure~\ref{fig_general} (a), (b), and (c) show the  trajectory data, the true coefficient for $u_{xx}$, and that for $uu_x$, respectively.

\textbf{GP-IDENT result:} 
 We use $4$ bases in space and $7$ bases in time to approximate the coefficients. Figure~\ref{fig_general} (d) shows the absolute error of the  trajectory simulated by the identified PDE by GP-IDENT, which is close to the true trajectory.
 Figure~\ref{fig_general} (e) and (f), display the relative $L_1$ errors ($\%$) for the reconstructed coefficients of $u_{xx}$ and $uu_x$, respectively. These figures demonstrates an accurate coefficient recovery of~\eqref{eq_exp1}. In (g)-(i), we show the absolute error of the simulated trajectory and the relative $L_1$ errors ($\%$) of the reconstructed coefficients when the given data have $2\%$ noise and SDD-9 is used for denoising. 
GP-IDENT successfully identified the underlying PDE, and the simulated trajectory remains close to the true one. We note that the reconstructed coefficient for $u_{xx}$ deviates from the true ones when $t\in (0.015,0.02)$ because the observed trajectory in (a) is mostly flat in this region. The flatness (derivatives being close to zero)  causes  a lack of local dynamics and leads to numerical instability.
The coefficient identification on this region is ill-posed. 

\begin{figure}[t!]
		\centering
		\begin{tabular}{cc}
			(a)&(b)\\
   \includegraphics[width = 0.4\textwidth,height=3.6cm]{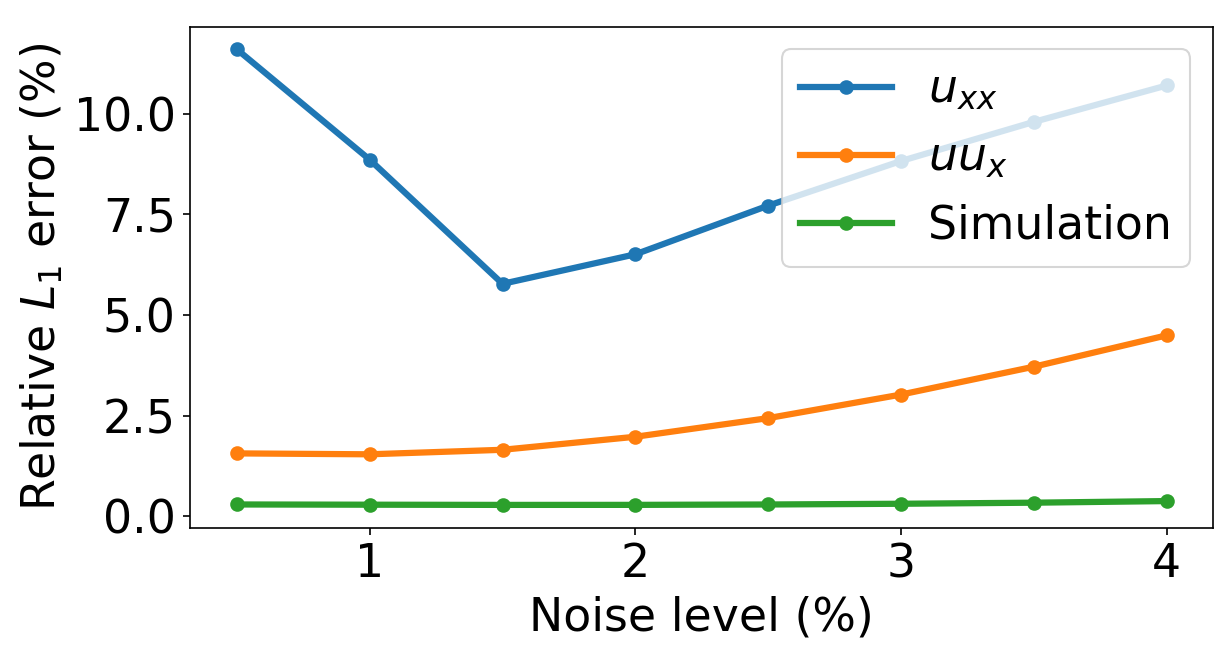}&
 \includegraphics[width = 0.45\textwidth,height=3.8cm]{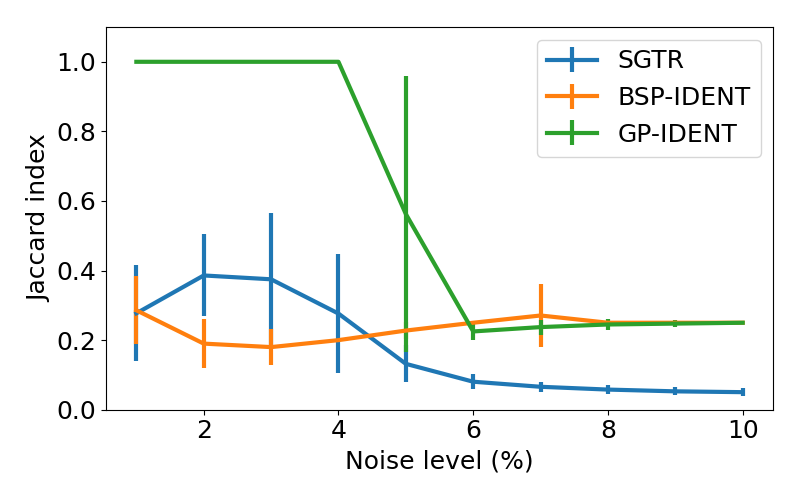}
		\end{tabular}
\caption{Varying noise level comparison for the viscous Burgers equation~\eqref{eq_exp1}: (a) Relative $L_1$ error for the coefficient for $u_{xx}$ (blue) and $u u_x$ (orange).  Green curve shows relative $L_1$ error between the true trajectory and the simulated trajectory of the identified model by GP-IDENT at various noise levels. Although the coefficient error for $u_{xx}$ is large due to ill-posedness, the simulated trajectory matches with less than 1\% error. 
(b) Accuracy of support identification measured by Jaccard index under various levels of noise. For each noise level, we ran $20$ independent experiments using the default dictionary with $56$ terms. For BSP-IDENT and GP-IDENT, we used SDD-9 for all levels of noise.}\label{fig_stats_Burgers}
	\end{figure}

\textbf{Robustness against various level of noise:} 
We demonstrate the robustness of GP-IDENT and compare with SGTR and BSP-IDENT for various noise levels. Figure~\ref{fig_stats_Burgers} (a) shows the relative $L_1$ coefficient error for $uu_x$ and $u_{xx}$, and the relative $L_1$ error between simulated trajectory using the model identified by GP-IDENT and the true trajectory (green). The coefficient identification  for $uu_x$ is robust to noise, yet
the coefficient reconstruction for $u_{xx}$ is an ill-posed problem, since the dynamics are flat at some regions, as shown in Figure \ref{fig_general}. Despite that the coefficient error for $u_{xx}$ is relatively large, the simulated trajectory matches the PDE solution with less than $1\%$ error. The reduction of the error in the early stage is caused by over-smoothing of SDD-9 when the noise level is very low.

Figure~\ref{fig_stats_Burgers} (b) shows the  Jaccard index between the exact support and the recovered one by SGTR, BSP-IDENT and GP-IDENT.  Both SGTR and BSP-IDENT fail to identify the correct features, while GP-IDENT successfully finds the correct model when the noise is below $4\%$.
  
 \subsection{Advection-diffusion equation with space-dependent coefficients}\label{sec_adv}
 
Consider the following advection-diffusion equation~\cite{Rudy2019DataDrivenIO} with spatially dependent coefficients, for $x\in[-5,5)$, and $t\in (0,5]$,
\begin{align}
u_t(x,t) &= \partial_x(a(x)u)+0.1u_{xx} = \partial_xa(x)u+a(x)u_x +0.1u_{xx} ~\label{eq_exp0}
	\end{align}
with initial condition $u(x,0) = \cos(2\pi x/5)$, and $a(x) = -1.5+\cos(2\pi x/5)$. This PDE is solved over a $256\times 256$ (space $\times$ time) grid. 

When the given data are noisy, SDD plays a critical role.  We show in Appendix \ref{Asec_sdd} Figure~\ref{fig_SDD} that noise is significantly amplified in the  finite difference scheme; whereas SDD effectively suppresses the perturbation in partial derivatives, thus it helps to identify the true dynamics. 

\begin{figure}[t!]
		\centering
		\begin{tabular}{ccc}
			(a)&(b)&(c)\\
   \includegraphics[width = 0.3\textwidth]{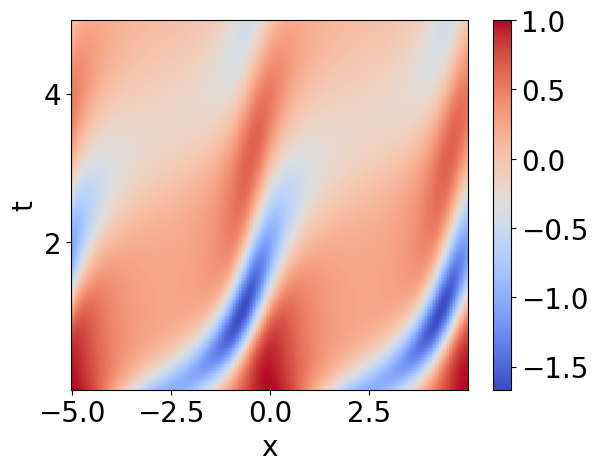}&\includegraphics[width = 0.3\textwidth]{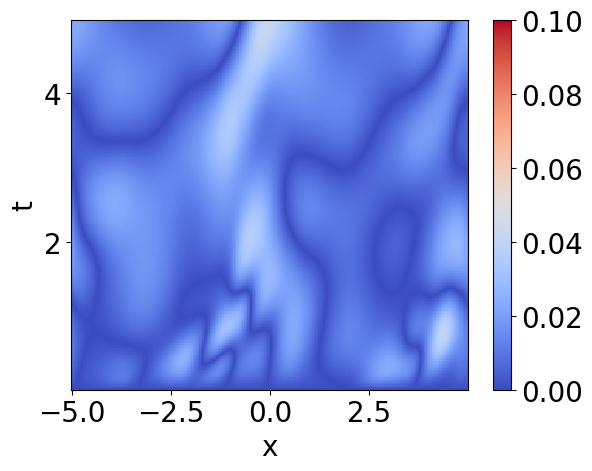}&\includegraphics[width = 0.3\textwidth]{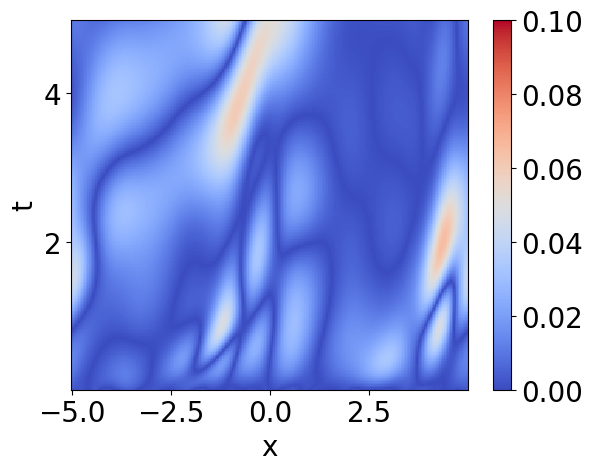}\\
   (d)&(e)&(f)\\
   \includegraphics[width=0.3\textwidth]{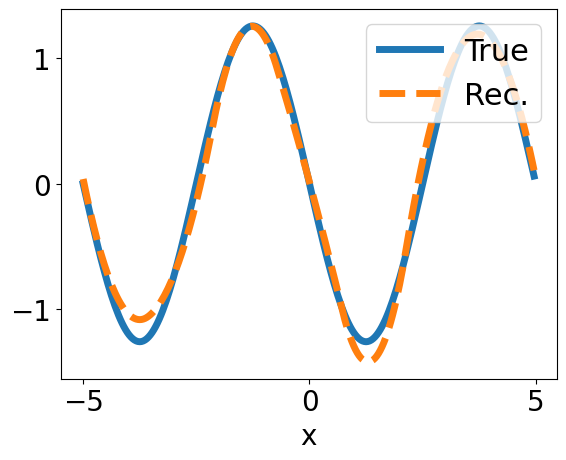}&
   \includegraphics[width=0.3\textwidth]{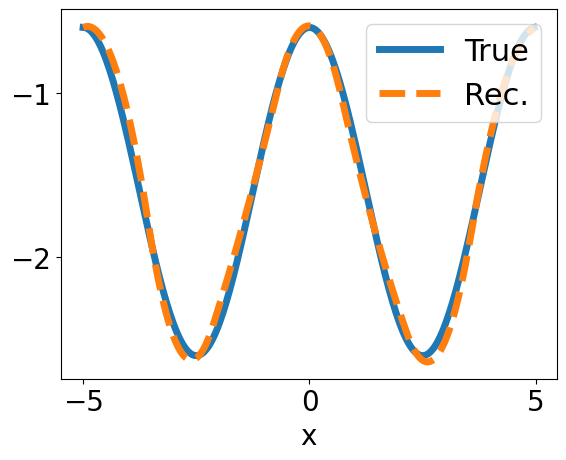}&
   \includegraphics[width = 0.3\textwidth]{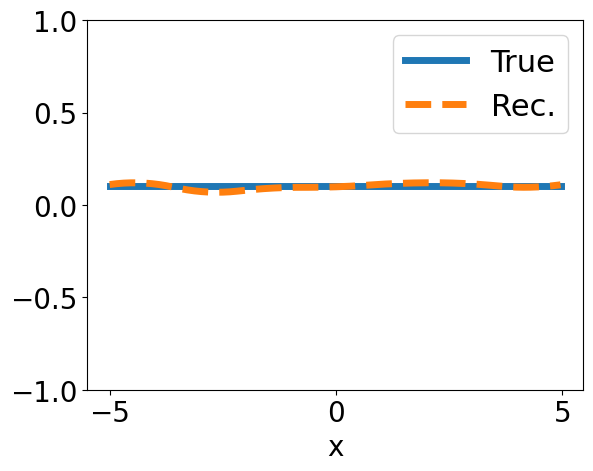}\\
   (g)&(h)&(i)\\
    \includegraphics[width=0.3\textwidth]{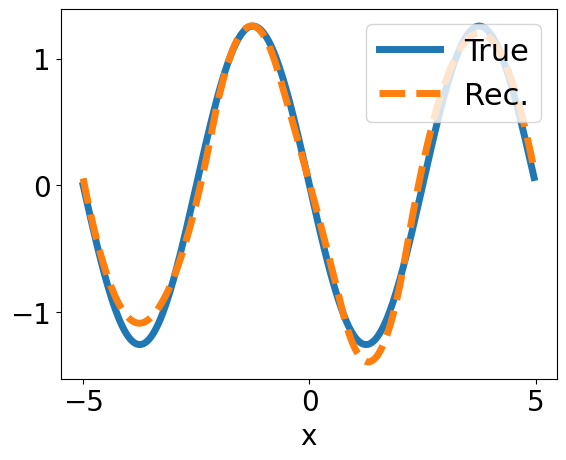}&
   \includegraphics[width=0.3\textwidth]{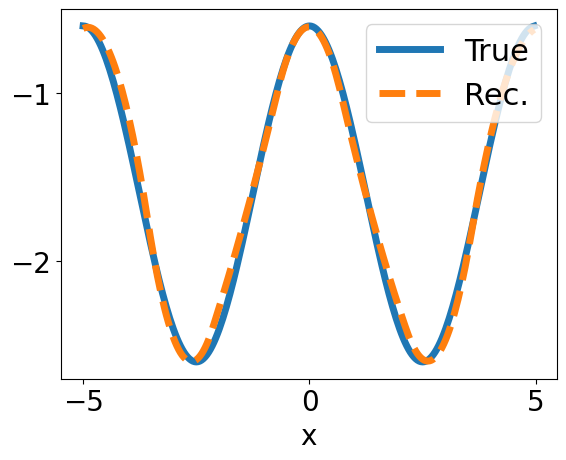}&
   \includegraphics[width = 0.3\textwidth]{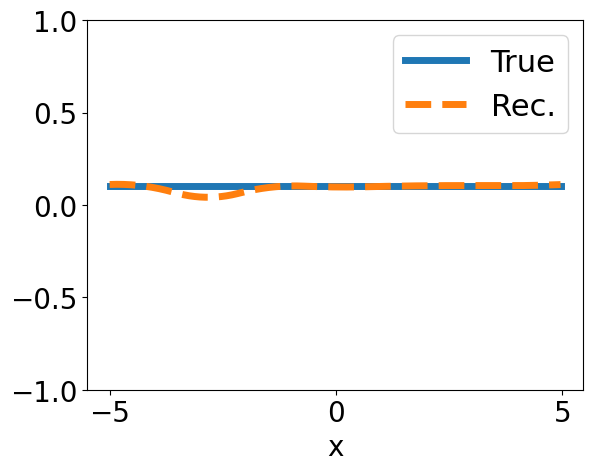}
		\end{tabular}
	\caption{GP-IDENT result for the advection-diffusion equation~\eqref{eq_exp0}: (a) observed clean trajectory. (b) and the second row (d)-(f) shows results from the clean data, and (c) and the third row (g)-(i) shows results from the given data with $1\%$ noise (SDD-15 is applied for denoising).   The first row shows absolute difference between the true (a) and the trajectory simulated by GP-IDENT. 
(d) and (g) are reconstruction of the coefficient of $u$, (e) and (h) of  $u_x$,  (f) and (i) of $u_{xx}$.  
}\label{fig_adv1}
\end{figure}

\textbf{GP-IDENT result:} For this experiment, we assume that we a priori know coefficients are only varying in space. Using $7$ bases in space for the coefficient approximation, GP-IDENT successfully identified the equation~\eqref{eq_exp0}.  Figure~\ref{fig_adv1} (b) and the second row, (d)-(f) show the reconstruction results with clean data, which stay close to the true coefficient values. 
We also test GP-IDENT when the data has $1\%$ noise using SDD-15 (Section~\ref{ss:sdd}).(c) shows the absolute error of the simulated trajectory, and the third row, (g)-(i) show the reconstructed coefficients. GP-IDENT shows robust recovery.  

\textbf{Stability against dictionary sizes:} 
We present the results with three  dictionaries of different sizes. In Dictionary I, we include partial derivatives of $u$ up to order $3$ and their products of no more than $3$ terms, in total of $35$ features. In Dictionary II, we include partial derivatives of $u$ up to order $4$ and their products of no more than $3$ terms, in total of $56$ features. In Dictionary III, we include partial derivatives of $u$ up to order $6$ and their products of no more than $4$ terms, in total of $330$ features. With each of these dictionaries, we apply GLASSO, SGTR, rSGTR, BSP-IDENT and GP-IDENT to identify~\eqref{eq_exp0} from a trajectory of data with or without noise. Table~\ref{tab_dic} compares the identified features of these methods with different noise levels. In this example, GLASSO does not converge when  Dictionary III is used, and except for this, all methods have correctly identified the true PDE, when the data has no noise.   When the given data have $1\%$ noise, GLASSO identifies the correct features for Dictionary II but not Dictionary I, which suggests that Dictionary II is more co-linear. SGTR fails to identify the correct terms in all cases. We note that in~\cite{Rudy2019DataDrivenIO},~\eqref{eq_exp0} is identified with a smaller dictionary. We find that rSGTR has identical performances as SGTR in terms of feature selection.  Both BSP-IDENT and GP-IDENT  yield the correct model.

\begin{table}[t!]
\centering
\begin{tabular}{c|c|c|c|c|c}
\hline
&\multicolumn{5}{c}{No noise}\\\hline
Method&GLASSO&SGTR& rSGTR&BSP-IDENT& GP-IDENT\\\hline
Dict. I & $\boldsymbol{u,u_x,u_{xx}}$&$\boldsymbol{u,u_x,u_{xx}}$&$\boldsymbol{u,u_x,u_{xx}}$&$\boldsymbol{u,u_x,u_{xx}}$& $\boldsymbol{u,u_x,u_{xx}}$\\\hline
Dict. II &$\boldsymbol{u,u_x,u_{xx}}$ &$\boldsymbol{u,u_x,u_{xx}}$&$\boldsymbol{u,u_x,u_{xx}}$&$\boldsymbol{u,u_x,u_{xx}}$& $\boldsymbol{u,u_x,u_{xx}}$\\\hline
Dict. III & $-$&$\boldsymbol{u,u_x,u_{xx}}$&$\boldsymbol{u,u_x,u_{xx}}$&$\boldsymbol{u,u_x,u_{xx}}$& $\boldsymbol{u,u_x,u_{xx}}$\\\hline
&\multicolumn{5}{c}{$1\%$ noise}\\\hline
Method&GLASSO&SGTR& rSGTR&BSP-IDENT& GP-IDENT\\\hline
Dict. I & 4 terms&$u,u_x$&$u,u_x$&$\boldsymbol{u,u_x,u_{xx}}$& $\boldsymbol{u,u_x,u_{xx}}$\\\hline
Dict. II &$\boldsymbol{u,u_x,u_{xx}}$ &$u,u_x$&$u,u_x$&$\boldsymbol{u,u_x,u_{xx}}$& $\boldsymbol{u,u_x,u_{xx}}$\\\hline
Dict. III & $-$&$u,u_x$&$u,u_x$&$\boldsymbol{u,u_x,u_{xx}}$& $\boldsymbol{u,u_x,u_{xx}}$\\\hline
&\multicolumn{5}{c}{$3\%$ noise}\\\hline
Method&GLASSO&SGTR& rSGT&BSP-IDENT& GP-IDENT\\\hline
Dict. I&$u,u_x,\partial_x^3u$& 5 terms&5 terms&$\boldsymbol{u,u_x,u_{xx}}$& $\boldsymbol{u,u_x,u_{xx}}$\\\hline
Dict. II &$u,u_x,\partial_x^3u$&5 terms&5 terms&$\boldsymbol{u,u_x,u_{xx}}$& $\boldsymbol{u,u_x,u_{xx}}$\\\hline
Dict. III & $-$& 5 terms& 5 terms&$\boldsymbol{u,u_x,u_{xx}}$& $\boldsymbol{u,u_x,u_{xx}}$\\\hline
&\multicolumn{5}{c}{$6\%$ noise}\\\hline
Method&GLASSO&SGTR& rSGTR&BSP-IDENT& GP-IDENT\\\hline
Dict. I& $u,u_x,\partial_x^3u$&18 terms&18 terms&$\boldsymbol{u,u_x,u_{xx}}$& $\boldsymbol{u,u_x,u_{xx}}$\\\hline
Dict. II &$u,u_x,\partial_x^3u$&10 terms&10 terms&$\boldsymbol{u,u_x,u_{xx}}$& $\boldsymbol{u,u_x,u_{xx}}$\\\hline
Dict. III&$-$&6 terms&6 terms&$\boldsymbol{u,u_x,u_{xx}}$& $\boldsymbol{u,u_x,u_{xx}}$\\\hline
\end{tabular}
\caption{The advection-diffusion equation~\eqref{eq_exp0} identification comparisons: GLASSO~\cite{yuan2006model}, SGTR~\cite{Rudy2019DataDrivenIO}, rSGTR~\cite{li2020robust}, BSP-IDENT, and GP-IDENT with three dictionaries and various noise levels. For  BSP-IDENT and GP-IDENT, SDD-15 is applied for denoising. Dictionary I has 35 features, II has 56 features, and  III has 330 features. Correct support identifications are marked in bold. GLASSO does not converge when Dictionary III is used.  GP-IDENT consistently identifies the correct terms. 
}\label{tab_dic}
\end{table}

\begin{figure}[t!]
		\centering
		\begin{tabular}{cc}
			(a)&(b)\\
\includegraphics[width = 0.45\textwidth]{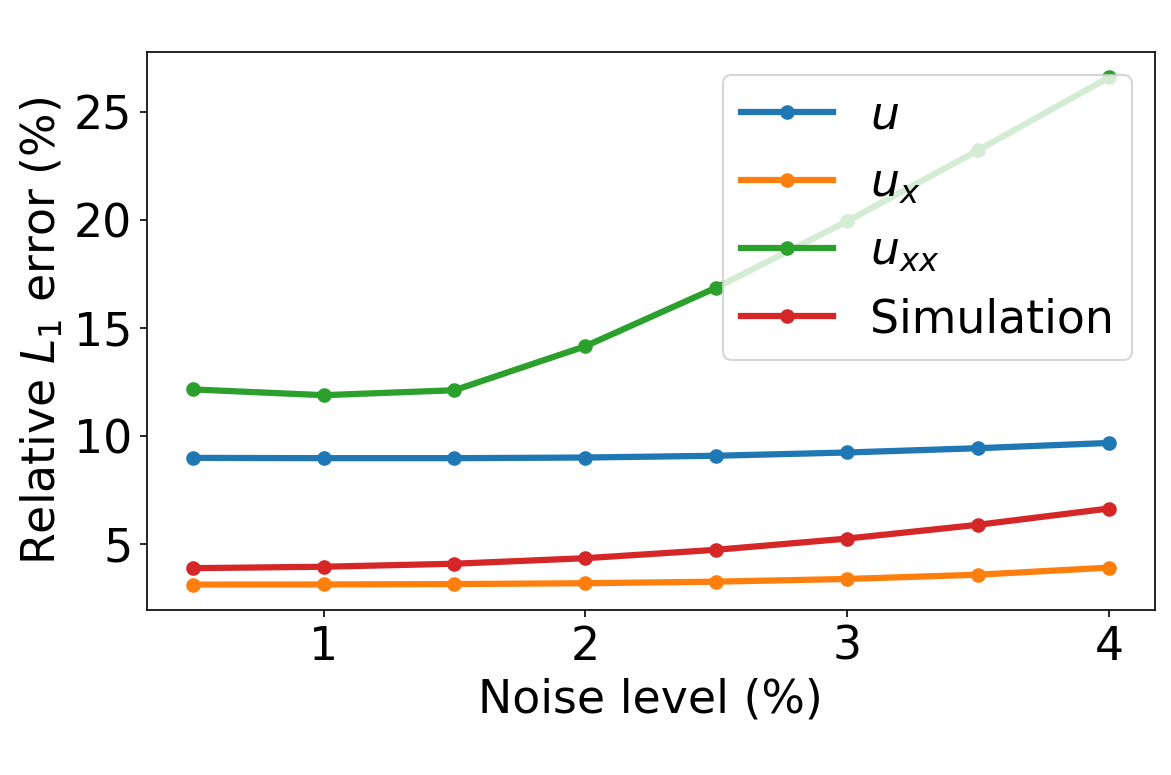}
 &  \includegraphics[width = 0.45\textwidth]{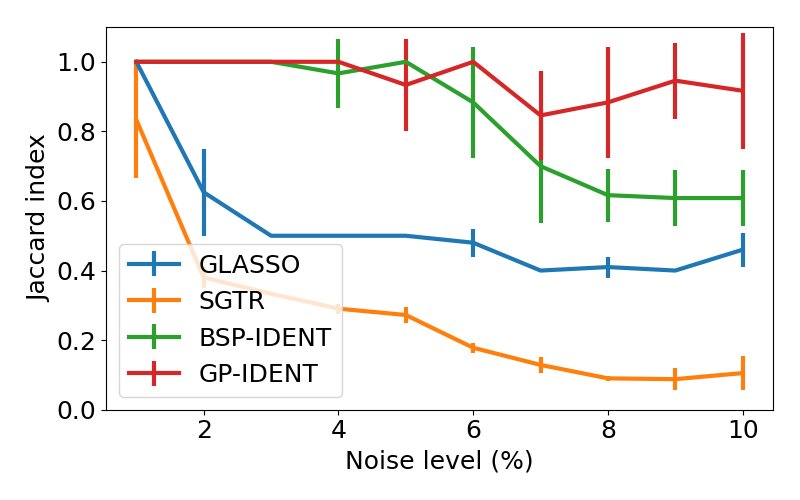}
		\end{tabular}
\caption{The advection-diffusion equation~\eqref{eq_exp0} identification comparisons with varying noise with Dictionary II.  (a) elative $L_1$ error for the coefficient for $u$ (blue), $u_x$ (orange), and $u_{xx}$ (green).  Red curve shows relative $L_1$ error between the true trajectory and 
the simulated trajectory of the identified model by GP-IDENT at various noise levels.  (b) Jaccard index showing correct support  identification with various levels of noise. For BSP-IDENT and GP-IDENT, we used SDD-15 at all levels of noise. 
 }\label{fig_stats_comp}
	\end{figure}

\textbf{Robustness against noise: } 
We demonstrate the robustness of GP-IDENT and compare with SGTR and BSP-IDENT for various noise levels. Figure~\ref{fig_stats_comp} (a) shows the relative $L_1$ coefficient error for $u$, $u_x$, and $u_{xx}$, and the relative $L_1$ error between simulated trajectory using the model identified by GP-IDENT and the true trajectory (red). The coefficient identification  for $u_x$ is robust to noise. Analogous to the case of Burgers' equation, the coefficient reconstruction for $u_{xx}$ is more challenging. Despite that the coefficient error for $u_{xx}$ is relatively large, the simulated trajectory closely matches the PDE solution.

Figure~\ref{fig_stats_comp} (b) shows the Jaccard indices of the identified features by these methods with various noise levels when Dictionary II is used. Overall GP-IDENT and BSP-IDENT  outperform the other methods. When the noise level is high, GP-IDENT yields better results than BSP-IDENT.

\textbf{Computational Efficiency:} Table  \ref{tab_eff} shows computational efficiency comparisons among the various methods for clean data.  For BSP-IDENT and GP-IDENT, left column of each method records the time when $K_{\max} = 10$ and the right column for $K_{\max} = 15$.  GP-IDENT and BSP-IDENT show fast converge.   We show the comparison between BSP and GPSP in Appendix \ref{Asec_bsp_Gpsp}.

\begin{table}[t!]
\centering
\begin{tabular}{c|c|c|c|c|c|c|c}
\hline
&\multicolumn{7}{c}{Identification time (sec)}\\\hline
Method&GLASSO&SGTR& rSGTR&\multicolumn{2}{c|}{BSP-IDENT}&\multicolumn{2}{c}{GP-IDENT}\\
\hline
$K_{\max}$& $-$& $-$&$-$ & 10 & 15 & 10 & 15 \\
\hline
Dict. I &$258.49$&$7.88$&$9.45$&$4.06$&$7.39$&$2.71$&$7.00$\\\hline
Dict. II&$356.34$&$12.77$&$15.19$&$2.98$&$6.98$&$3.61$& $9.24$ \\\hline
Dict. III&$-$&$145.60$&$229.46$&$5.73$&$9.82$&$4.19$&$8.30$\\\hline
\end{tabular}
\caption{Comparison of computational efficiency for the advection-diffusion equation \eqref{eq_exp0} among GLASSO~\cite{yuan2006model}, SGTR~\cite{Rudy2019DataDrivenIO}, rSGTR~\cite{li2020robust}, BSP-IDENT, and GP-IDENT for clean data and different dictionary sizes. GLASSO fails to converge when Dictionary III is used.  GP-IDENT and BSP-IDENT show fast converge.}\label{tab_eff}
\end{table}

\subsection{Fisher's equation with time-dependent coefficients}
\begin{figure}
		\centering
		\begin{tabular}{ccc}
			(a)&(b)&(c)\\
			\includegraphics[width = 0.3\textwidth]{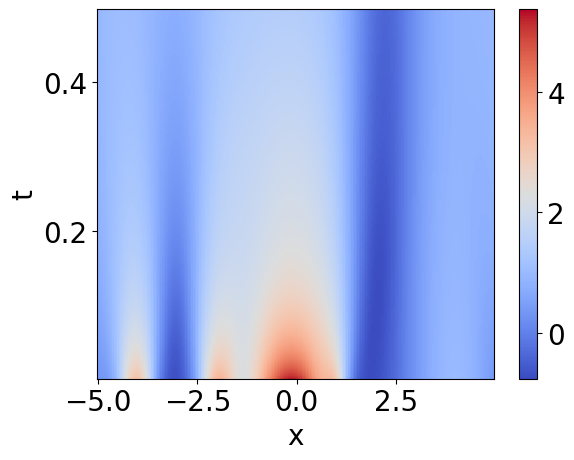}&\includegraphics[width = 0.32\textwidth]{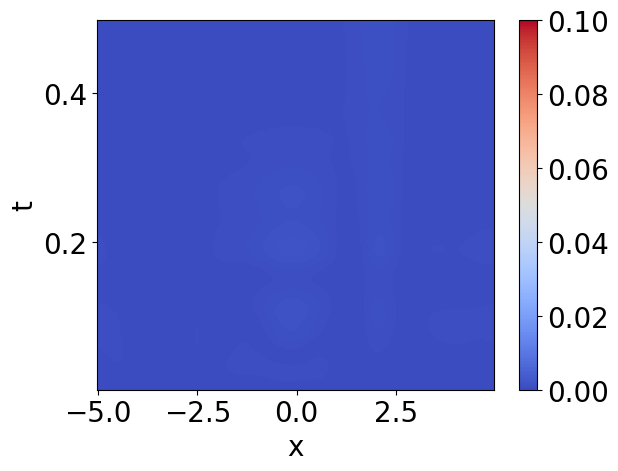}&\includegraphics[width = 0.32\textwidth]{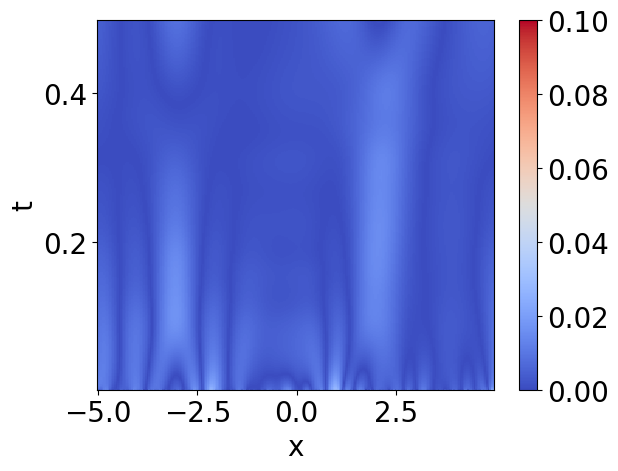}\\
   (d)&(e)&(f)\\
   \includegraphics[width=0.3\textwidth]{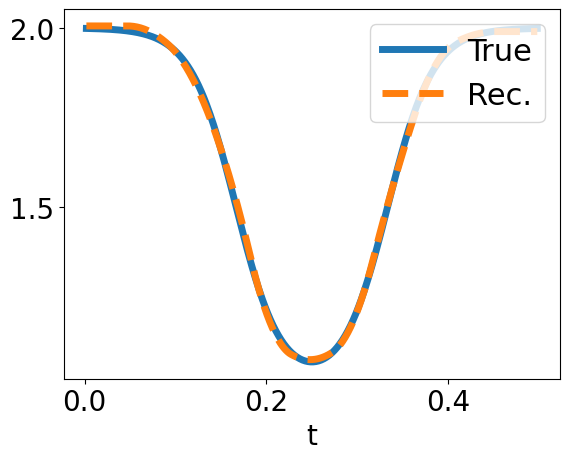}&
   \includegraphics[width=0.3\textwidth]{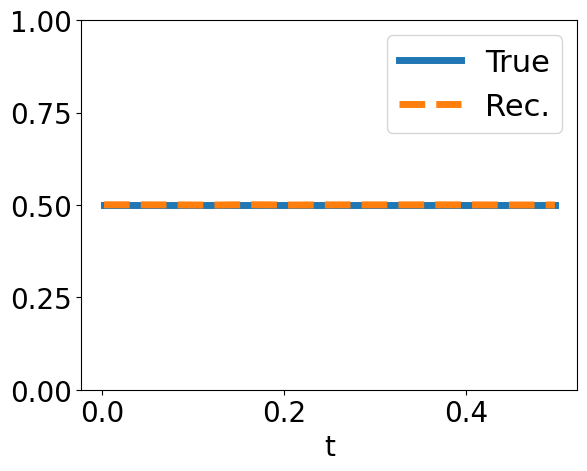}&
   \includegraphics[width = 0.3\textwidth]{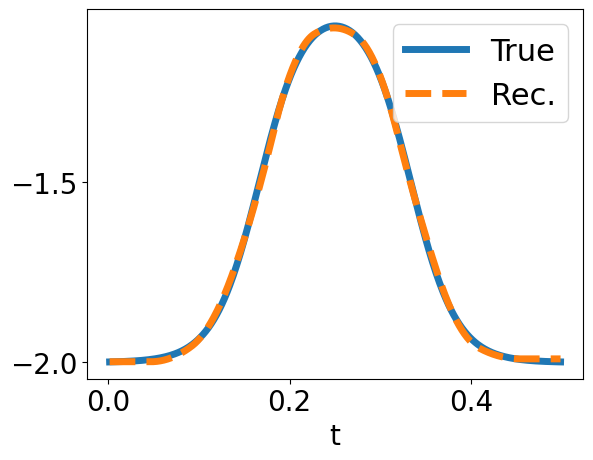}\\
   (g)&(h)&(i)\\
   \includegraphics[width=0.3\textwidth]{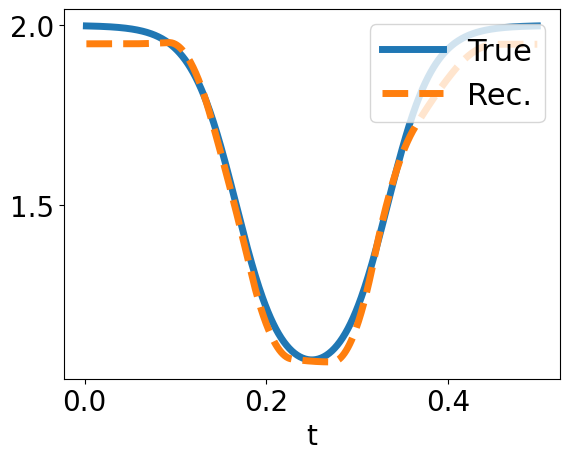}
   &\includegraphics[width=0.3\textwidth]{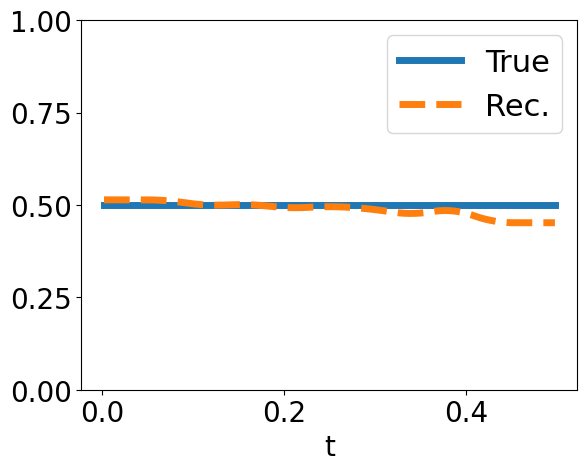}
   &\includegraphics[width = 0.3\textwidth]{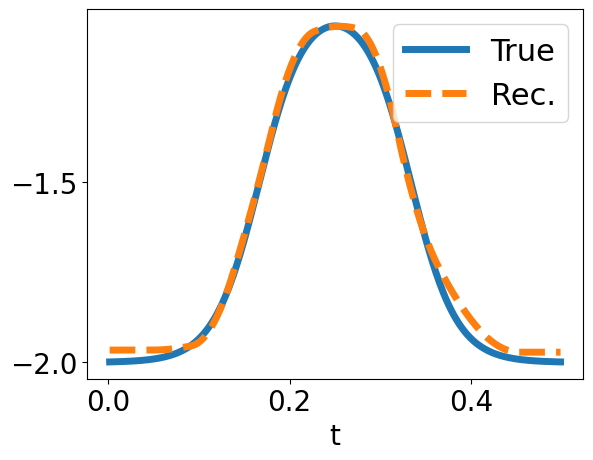}
		\end{tabular}
\caption{GP-IDENT result on the Fisher's equation~\eqref{eq_exp_fisch}: 
(a) observed clean trajectory. (b) and the second row (d)-(f) shows GP-IDENT results from the clean data, and (c) and the third row (g)-(i) shows results from the given data with $2\%$ noise (SDD-15 is applied for denoising).   The first row shows absolute difference between the true (a) and the trajectory simulated by GP-IDENT.  (d) and (g) are reconstruction of the coefficient of $u$, (e) and (h) of  $u_{xx}$,  (f) and (i) of $u^2$.}\label{fig_fisch1}
\end{figure}
Consider the Fisher's equation with time-dependent growth rate~\cite{ougun2007exact} widely studied in physics and genetics
\begin{align}
		u_t(x,t) &= 0.5u_{xx}(x,t) +a(t)u(x,t)(1-u(x,t))\;,~x\in[-5,5), t\in (0,0.8]\label{eq_exp_fisch}
\end{align}
where 
\begin{align}
a(t) = 1 + \tau_{-}\left(t;s,\frac{0.8}{3}\right)+\tau_{+}\left(t;s,\frac{1.6}{3}\right)
\end{align}
and
\begin{align}
		\tau_{\pm}(t;s,t_b) = \frac{1}{2}+\frac{1}{2}\tanh\left(\pm \frac{s(t-t_b)}{T_{\max}}\right),~t\in [0,T_{\max}]\label{eq_trans}
	\end{align}
which reflects a smooth transition with rate $s$ between different states separated by the break point $t_b$. We take the initial condition
\begin{align}
u(x,0)=5e^{-x^2}+3e^{-(2x+4)^2}+2e^{-(3x-3)^2}+4e^{-(2x+8)^2}+\cos(4(x+1)\pi/10),
\end{align}
and numerically solve it on a $256\times 512$ grid for  $s= 10$.
 
\textbf{GP-IDENT Result:} For this experiment, we assume that we a priori know coefficients are only varying in time. We apply GP-IDENT with $9$ bases in time to approximate the coefficients. Figure~\ref{fig_fisch1} (a) shows the clean trajectory, (d)-(f) present the identified coefficients compared to the true coefficients when the given data are clean, and (b) shows the absolute error of the trajectory simulated from the identified model. GP-IDENT identifies varying coefficients accurately. When the data have $2\%$ noise, we apply SDD-15 for denoising, and GP-IDENT  identifies the correct model. (c) shows the absolute error of the trajectory simulated from the identified model, and (g)-(i) display the identified coefficients. GP-IDENT is robust to noise. 

\textbf{Comparisons:} In Table~\ref{tab_fisch}, we compare GLASSO, SGTR, rSGTR, BSP-IDENT, and GP-IDENT for identifying Fisher's equation~\eqref{eq_exp_fisch} with the  default dictionary using clean and noisy data. When the given data are clean, all methods identify the correct model. For the data perturbed by noise, GLASSO identifies extra terms, and both SGTR and rSGRT fail to find the correct terms.  BSP-IDENT and GP-IDENT identify the correct model up to $3\%$ noise.  
\begin{table}[ht]
\centering
\begin{tabular}{c|c|c|c|c|c}
\hline
Method&GLASSO&SGTR& rSGTR&BSP-IDENT& GP-IDENT\\\hline
No noise& $\boldsymbol{u,u_{xx},u^2}$&$\boldsymbol{u,u_{xx},u^2}$&$\boldsymbol{u,u_{xx},u^2}$&$\boldsymbol{u,u_{xx},u^2}$& $\boldsymbol{u,u_{xx},u^2}$\\\hline
$1\%$ noise&4 terms&$u_{xx},u^2,u^3$&$u_{xx},u^2,u^3$&$\boldsymbol{u,u_{xx},u^2}$& $\boldsymbol{u,u_{xx},u^2}$\\\hline
$2\%$ noise&5 terms&6 terms&6 terms&$\boldsymbol{u,u_{xx},u^2}$& $\boldsymbol{u,u_{xx},u^2}$\\\hline
$3\%$ noise&5 terms&$u^2$&$u^2$&$\boldsymbol{u,u_{xx},u^2}$& $\boldsymbol{u,u_{xx},u^2}$\\\hline
\end{tabular}
\caption{Comparison result for Fisher's equation~\eqref{eq_exp_fisch}. Table shows identified features of GLASSO, SGTR, rSGTR, BSP-IDENT, and GP-IDENT for data with several levels of noise. For both BSP-IDENT and GP-IDENT, SDD-15 was applied for all levels of noise. Correct identifications are marked in bold. Both BSP-IDENT and GP-IDENT consistently identifies the correct terms. }\label{tab_fisch}
\end{table}

\section{Conclusion}\label{sec_conclusion}

We propose an effective and efficient  method, GP-IDENT, for identifying parametric PDEs with space and time-dependent coefficients. Our method generates a few candidates by a greedy algorithm called GPSP at various levels of group sparsity.
GPSP algorithm find a sparse solution to the feature system for any given group sparsity. After generating the candidates, we find the optimal model by considering the smallest sparsity $k$ for a small RR, i.e. $s_k < \rho$.  This motivates to find simple equations where RR does not reduce further by adding more complex terms. We demonstate the effectiveness and efficiency of GP-IDENT on various types of PDEs and compare it with the state-of-the-art methods for PDE identification with varying coefficients. In our experiments, GP-IDENT consistently yields accurate and robust results.

To further improve the identification accuracy especially under influence of high levels of noise, model selection criterion adaptive to noise level may need to be considered in the future. 
\bibliographystyle{abbrv}
\bibliography{bibfile}

\appendix

\section{Effect of SDD, and the experiment settings for Table \ref{tab_more_exp} space and time varying equations}\label{Asec_sdd}

\begin{figure}[t!]
		\centering
		\begin{tabular}{ccc}
			(a)&(b)&(c)\\
			\includegraphics[width = 0.3\textwidth]{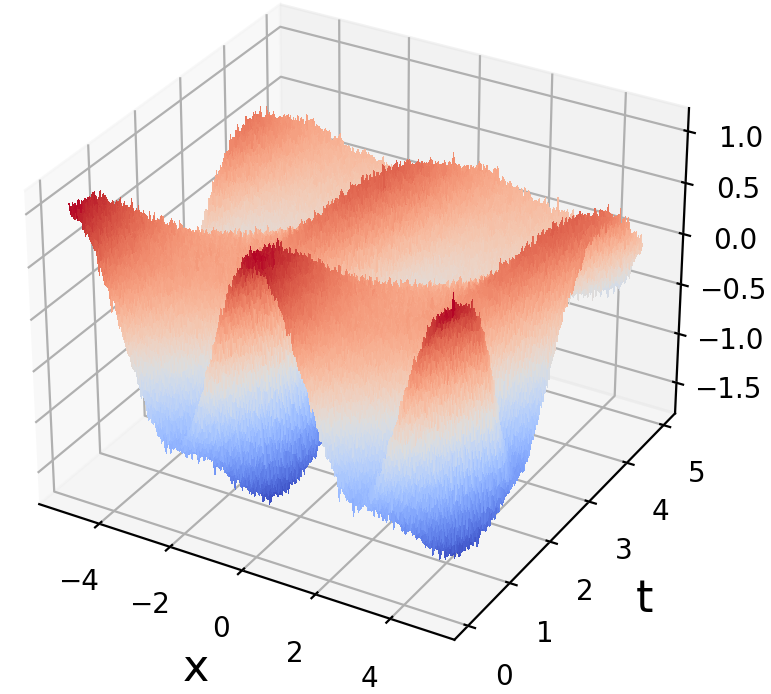}&
			\includegraphics[width=0.3\textwidth]{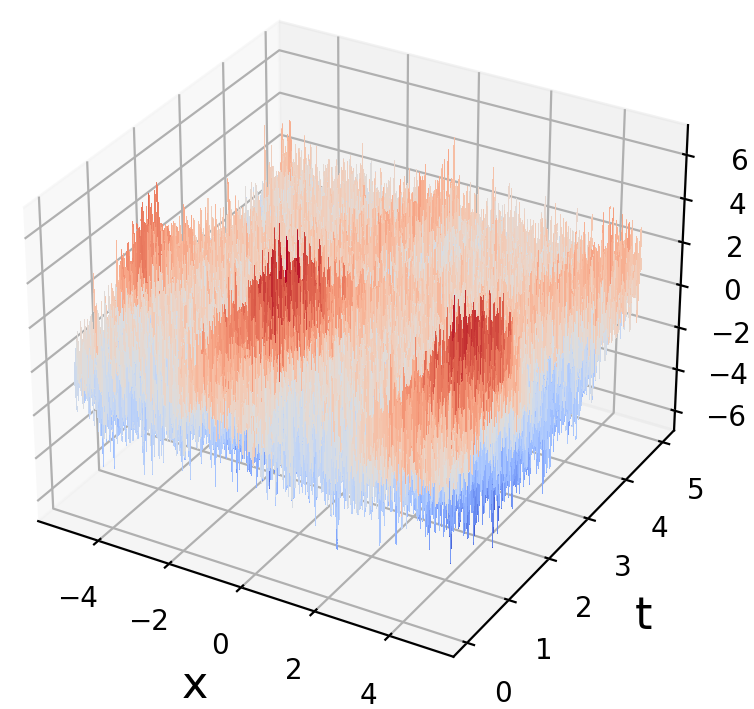}
   &
			\includegraphics[width=0.3\textwidth]{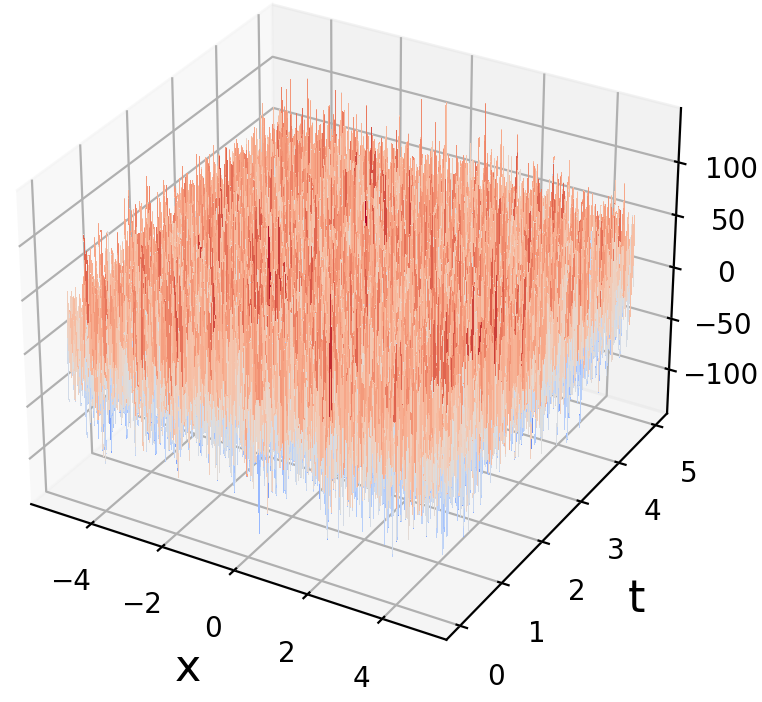}\\
   (d)&(e)&(f)\\
			\includegraphics[width = 0.3\textwidth]{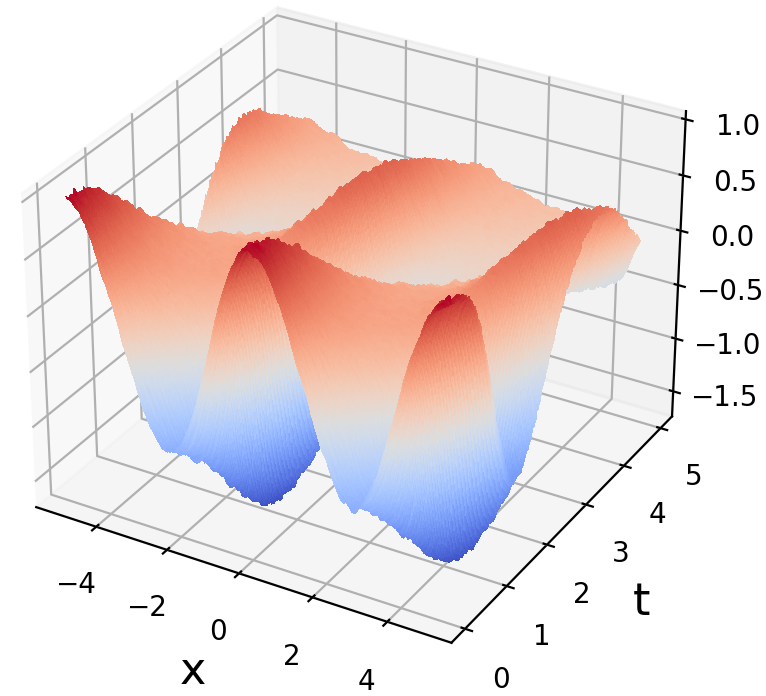}&
			\includegraphics[width=0.3\textwidth]{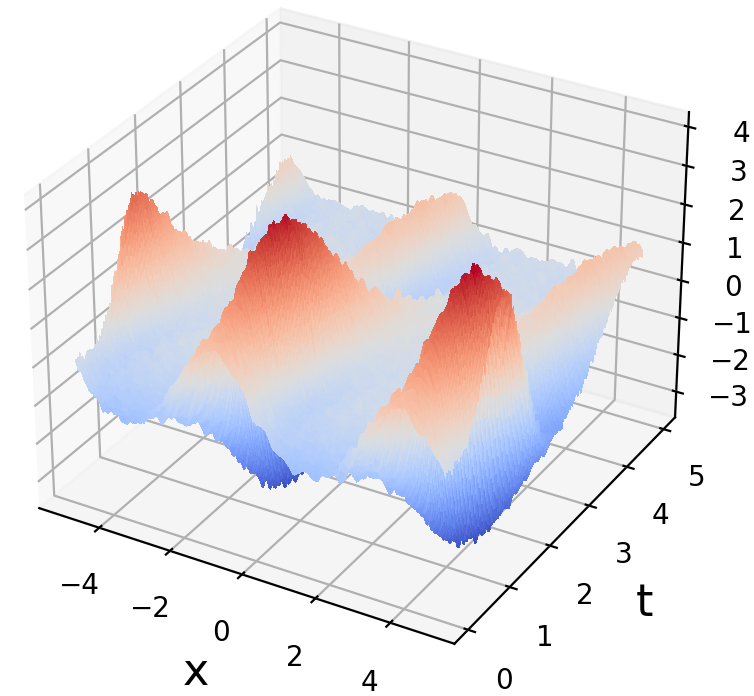}
   &
			\includegraphics[width=0.3\textwidth]{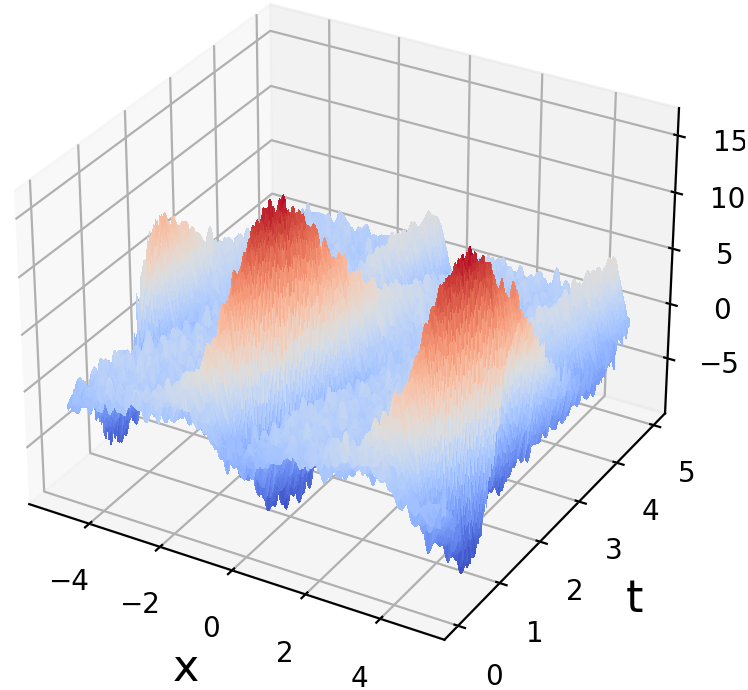}
		\end{tabular}
\caption{For advection-diffusion equation~\eqref{eq_exp0}, influence of noise and effectiveness of SDD: (a) A noisy trajectory with $10\%$ noise, (b) $u_{x}$ and (c)  $u_{xx}$ computed from the noisy data. With SDD-15 in the second row, (d) denoised $u$, (e) denoised $u_x$, and (f) denoised $u_{xx}$ are more stabilized. }\label{fig_SDD}
	\end{figure}

Figure \ref{fig_SDD} shows that when the given data are noisy, it is significantly amplified in the  finite difference scheme.  SDD effectively suppresses the perturbation in partial derivatives, and helps to identify the true dynamics. 


In Table \ref{tab_more_exp_info}, we present the details of these experiment settings including the coefficients, grid, number of bases, and window size for SDD for the experiments in Table \ref{tab_more_exp}. 

\begin{table}[t!]
\small
\centering
\begin{tabular}{c|c|c|c|c}
\hline
Model& Coefficients&Grid&Bases&SDD\\\hline
KdV&$\begin{cases}a(x) &= 0.5\cdot(2+0.3\cos(\pi x/2))\cdot\\
&(1+\tau_{+}(t;10,0.05))\\
b(x,t) &=0.01\cdot(0.5+0.1\sin(\pi x/2))\cdot\\
&(1+\tau_{-}(t;10,0.05))
\end{cases}$ & \parbox{2.5cm}{\begin{align*}&256\times512\\ &[-2,2)\times[0,0.1]\end{align*}}&$5,5$&$5$\\\hline
KS&$\begin{cases}a(x) &= 2+\sin(2\pi x/30)/4\\
b(x,t) &=(-1+e^{-(x-2)^2/5}/4)\cdot\\
&(2+\tau_{+}(t, 5, 30))\\
c(x,t)&=(-1-e^{-(x+2)^2/5}/4)\cdot\\
&(2+\tau_{+}(t, 5, 30))
\end{cases}$ & \parbox{2.5cm}{\begin{align*}&512\times512\\ &[-30,30)\times[0,60]\end{align*}}&$9,5$&$15$\\\hline
Sch&\parbox{4cm}{\begin{align*}a(x,t) &=-5\cos(\pi x/2)\cdot\\
&(0.5+\tau_{+}(t;5,0.2))\end{align*}} & \parbox{2.5cm}{\begin{align*}&100\times2000\\ &[-2,2)\times[0,2]\end{align*}} &$5,5$&$7$\\\hline
NLS&\parbox{5cm}{\begin{align*}a(x,t) &=(1+0.2\cos(\pi x/2))\cdot\\&(1+0.5\tau_{+}(t;5,0.2))\end{align*}} & \parbox{2.5cm}{\begin{align*}&100\times2000\\ &[-2,2)\times[0,0.5]\end{align*}}&$5,5$&$7$\\
\hline
\end{tabular}
\caption{Details of equations tested in Table~\ref{tab_more_exp}.  
The Grid column shows the space mesh size $\times $ the time mesh size on top, and the space domain $\times$ the time domain on the bottom. The Bases column  shows the number of bases used for space and time respectively. The SDD column  records the smoothing window size of SDD for the noisy data in each case. 
}\label{tab_more_exp_info}
\end{table}

\section{Effects of Reduction in Residual (RR) }\label{Asec_RR}

\begin{figure}[t!]
		\centering
  \begin{tabular}{cc}
  (a) & (b) \\
  \includegraphics[width = 0.44\textwidth]{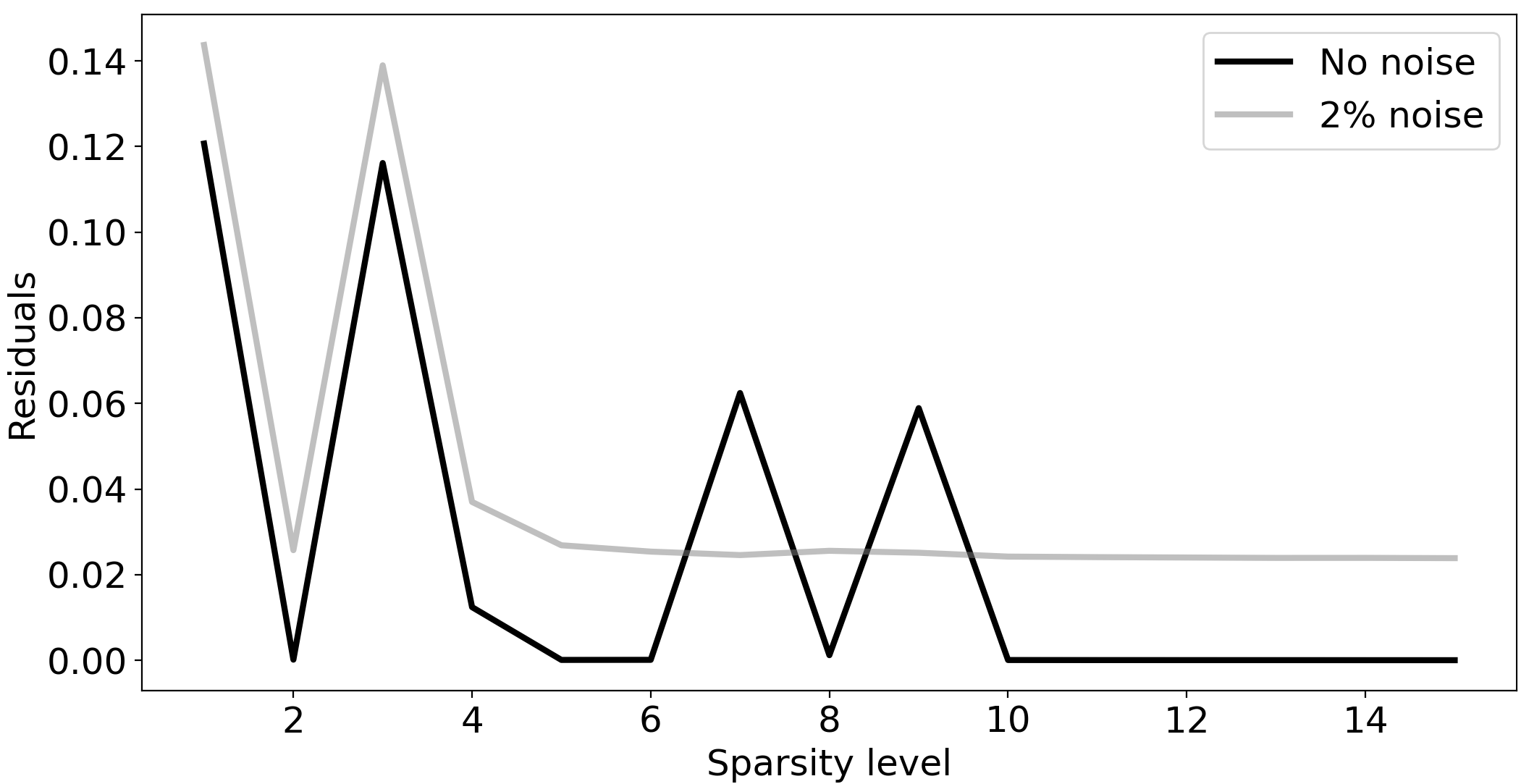} &
			\includegraphics[width = 0.47\textwidth]{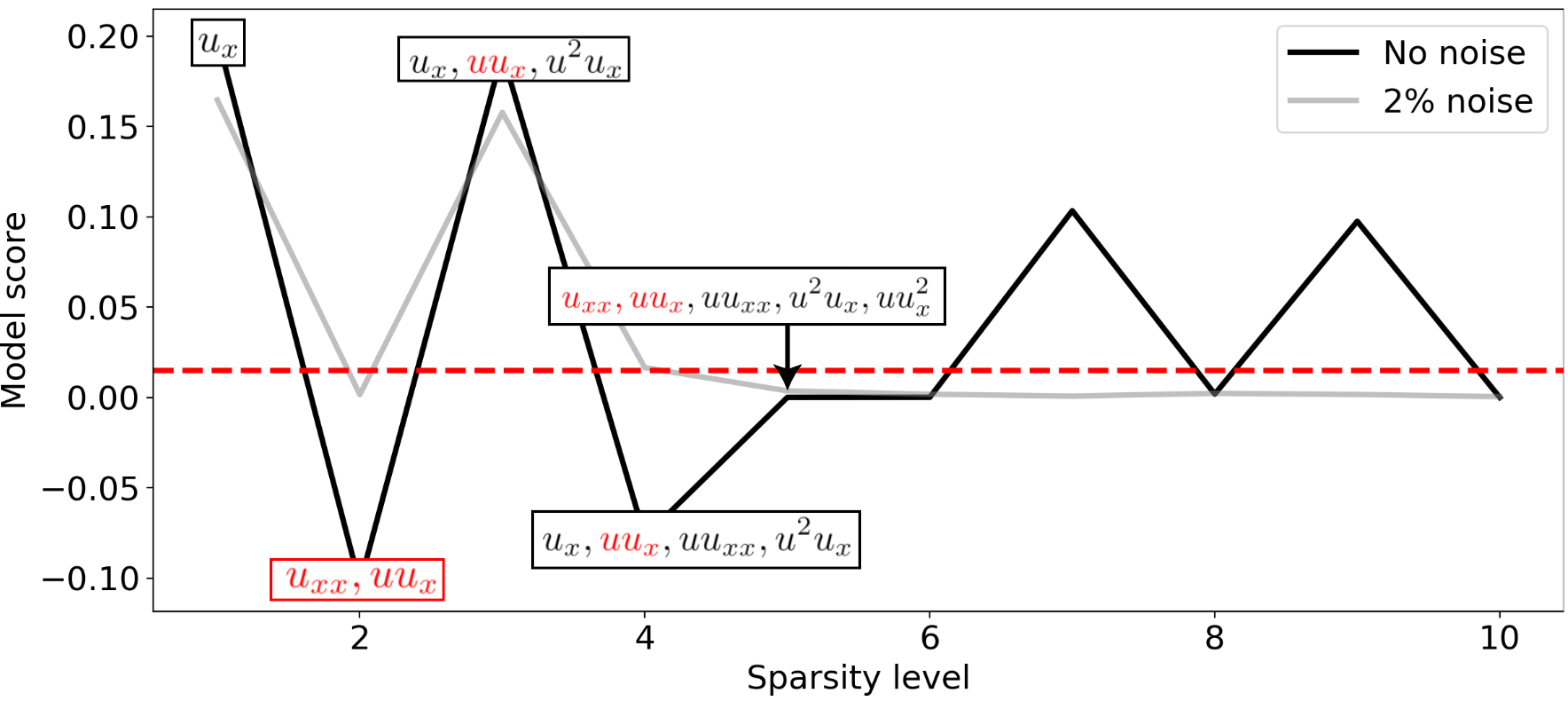}
   \end{tabular}
\caption{Effects of RR for the viscous Burgers equation~\eqref{eq_exp1}. For both graphs, the black curve is when there is no noise, and the gray curve is for 2\% noise. (a) Residuals of the candidate models from GPSP of various sparsity levels ($K_{\max}=15$). (b) RR score in ~\eqref{eq:cand_score} for candidates generated by GPSP using 
$L=5$. The red dashed curve represents the default threshold $\rho=0.015$, and the identified model is the one whose score first hits below $\rho$. The correct features are marked in red  in (b). 
}\label{fig_general2}
\end{figure}

In this paper, we propose the RR scores~\eqref{eq:cand_score} to select the identified PDE from a pool of candidates given by GPSP at various levels of sparsity. Using the Burgers equation~\eqref{eq_exp1} as an example, Figure~\ref{fig_general2} demonstrates that the RR scores are effective in selecting the correct model. 
For both (a) and (b), the black curves are when there is no noise, and the gray curves are for 2\% noise.  (a) shows residuals for each sparsity level.  As the sparsity level gets bigger, the residual curves fluctuate since different sparsity levels are produced by GPSP individually. For example, when the sparsity is $1$, the candidate contains $u_x$. When the sparsity is $2$, $u_x$ is removed, and the correct features $u_{xx}, uu_x$ are included. If the sparsity level is $3$, the correct feature $u_{xx}$ is not selected, which leads to an increment of residuals. The residual curves do not give clear indications about the optimal models. In (b), we show the RR curves as well as the threshold $\rho=0.015$ used in~\eqref{eq_thresh}  marked by the dashed red line. By our criterion, the optimal model matches the true one, since after the model with sparsity $2$, models with more complexity do not contribute to significant reduction in the residuals. Moreover, we note that when the given data have noise, the candidates' RR scores become less oscillatory, and the score for the correct model approaches the threshold $\rho=0.015$. This is commonly observed  in other PDEs as well. 

\section{Explanation of GPSP over BSP in PDE identification}\label{sec_GPSP_BSP_PDE}

Consider the transport equation with a constant speed $a\neq 0$
\begin{align}
u_t(x,t) = au_x(x,t)\label{eq:amb_eq1}
\end{align}
and its solution $f(x+at)$ for some smooth function $f$, which is nowhere zero.  If the hypothesis space contains $f(x+at)$ and the dictionary contains  $u_x, uu_x$, it is possible to confuse~\eqref{eq:amb_eq1} with
\begin{align}
u_t(x,t) = \frac{a}{f(x+at)}u(x,t)u_x(x,t)\;~\text{or}~u_t(x,t) = a_1u_x(x,t)+\frac{a_2}{f(x+at)}u(x,t)u_x(x,t)\label{eq:amb_eq2}
\end{align}
where $a_1 + a_2 = a$ and $a_1,a_2\neq 0$, in which case, all these PDE models are valid.

In practice, the dimension of the hypothesis space $M$ is finite, and the hypothesis space is confined by the resolution of the sampling grid for numerical stability.   Ideally, the PDE model with the least coefficient approximation error by the hypothesis space should be selected, and this is where GPSP differs from BSP.

We denote $g(x,t)=\frac{1}{f(x+at)}$ and decompose $g(x,t)=g_M(x,t)+e_M(x,t)$ where $g_M$ is the orthogonal projection of $g$ to $\mathcal{H}_M$, and $e_N\perp\mathcal{H}_M$ denotes the residual. For simplicity, we assume normalization is applied and  the dictionary is simply $\{u_x,uu_x\}$. We compare GPSP with BSP when the sparsity level is fixed at $1$, that is, each method selects just one feature, and we focus on the selection in the initial step. In BSP, we are  comparing
\begin{align}
 \sqrt{\sum_{m=1}^M\left(\langle u_t, B_mu_x\rangle\right)^2}\;~\text{with}~\sqrt{\sum_{m=1}^M\left(\langle u_t, B_muu_x\rangle\right)^2}
\end{align}
where $\{B_m\}_{m=1}^M$ are basis functions and the inner product is understood as operations over the grid points, for example,
\begin{align}
\langle u_t, B_mu_x\rangle=\sum_{i=1}^I\sum_{n=1}^Nu_t(x_i,t_n)B_m(x_i,t_n)u_x(x_i,t_n)
\end{align}
By the hypothesis space approximation, (\ref{eq:amb_eq1}) and (\ref{eq:amb_eq2}), we have
\begin{align}
\langle u_t, B_mu_x\rangle &= a\langle(g_M+e_M)uu_x,B_mu_x\rangle=a\langle g_M, B_muu_x^2\rangle+a\langle e_M, B_muu_x^2\rangle\\
\langle u_t, B_muu_x\rangle &= a\langle u_x,B_muu_x\rangle=a\langle g_M,B_muu^2_x\rangle+a\langle 1-g_M,B_muu^2_x\rangle
\end{align}
thus
\begin{align}
\langle u_t, B_muu_x\rangle-\langle u_t, B_mu_x\rangle=a\langle 1-g,B_muu^2_x\rangle.
\end{align}
It indicates that in the first step of BSP, the choice between $u_x$ and $uu_x$ is \textbf{independent} of the approximation error $e_M$; instead, the sign of $a$ as well as the magnitude of the trajectory affect the choice.  As for GPSP, we compare 
\begin{align}
\frac{\langle\text{Proj}(u_t, \text{span}_mB_mu_x),u_t\rangle}{\|\langle\text{Proj}(u_t, \text{span}_mB_mu_x)\|_2}=\frac{\langle\text{Proj}(au_x, \text{span}_mB_mu_x),u_t\rangle}{\|\langle\text{Proj}(au_x, \text{span}_mB_mu_x)\|_2}=\|u_t\|_2
\end{align}
with 
\begin{align}
\frac{\langle\text{Proj}(u_t, \text{span}_mB_muu_x),u_t\rangle}{\|\langle\text{Proj}(u_t, \text{span}_mB_muu_x)\|_2}\leq \|u_t\|_2\label{eq_proj_inn}
\end{align}
where $\text{Proj}(u_t, \text{span}_mB_mu_x)$ denotes the projection of $u_t$ to the column space spanned by $\{B_mu_x\}_{m=1}^M$. We note that in~\eqref{eq_proj_inn}, the equality holds if and only if $e_M=0$. Therefore, we conclude that GPSP will choose $u_x$ over $uu_x$ if the approximation error for the finite dimensional hypothesis space is non-zero.  In other words, the choice of GPSP is \textbf{dependent} on the approximation error.

\begin{figure}[t!]
		\centering
		\begin{tabular}{ccc}
			(a)&(b)&(c)\\
   \includegraphics[width = 0.3\textwidth]{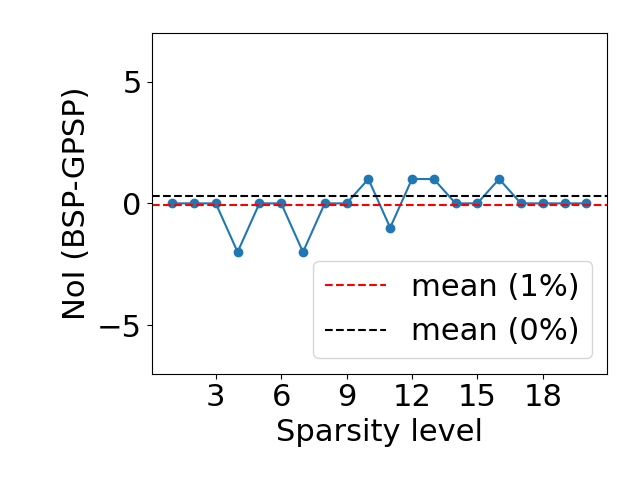}&
			\includegraphics[width=0.3\textwidth]{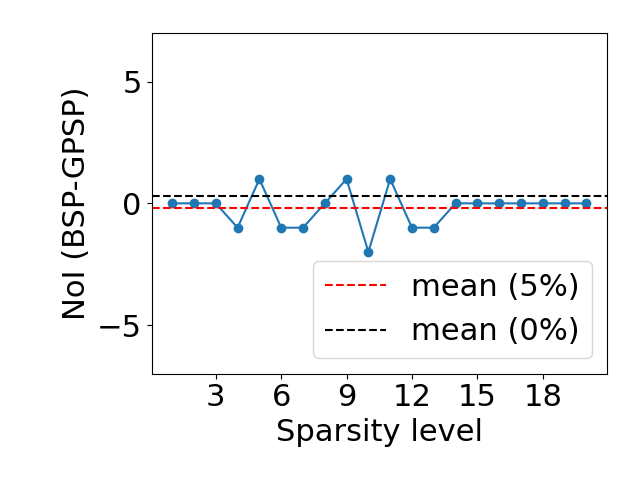}
   &
			\includegraphics[width=0.3\textwidth]{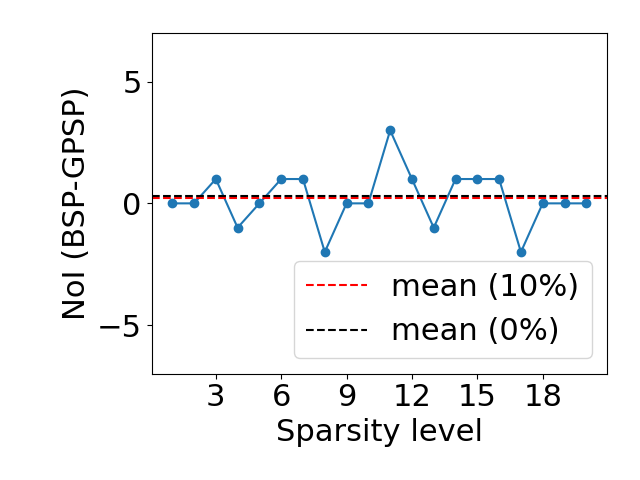}\\
                (d)&(e)&(f)\\
			\includegraphics[width = 0.3\textwidth]{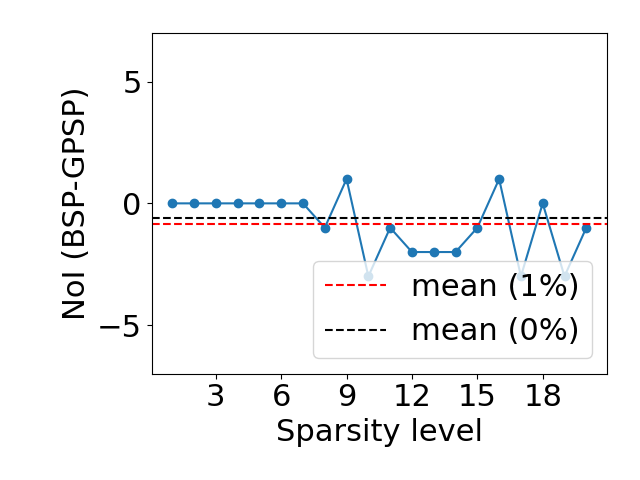}&
			\includegraphics[width=0.3\textwidth]{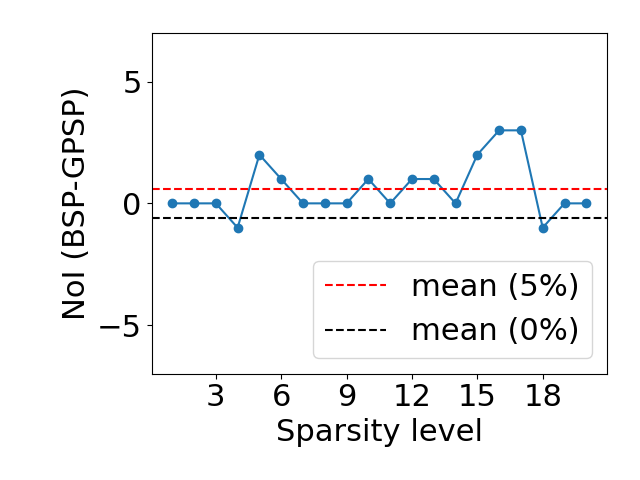}
   &
			\includegraphics[width=0.3\textwidth]{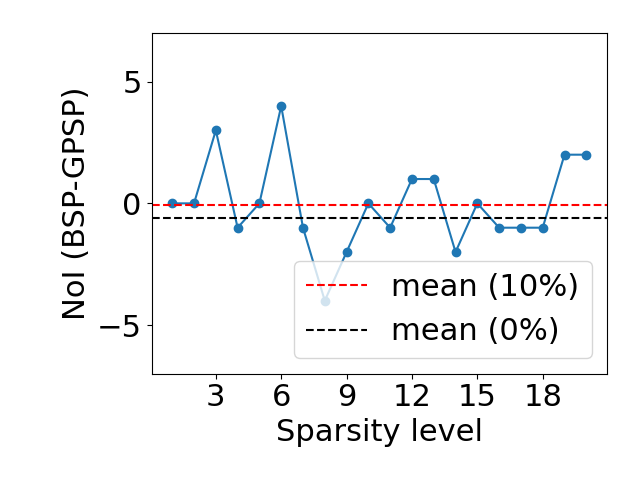}\\(g)&(h)&(i)\\
			\includegraphics[width = 0.3\textwidth]{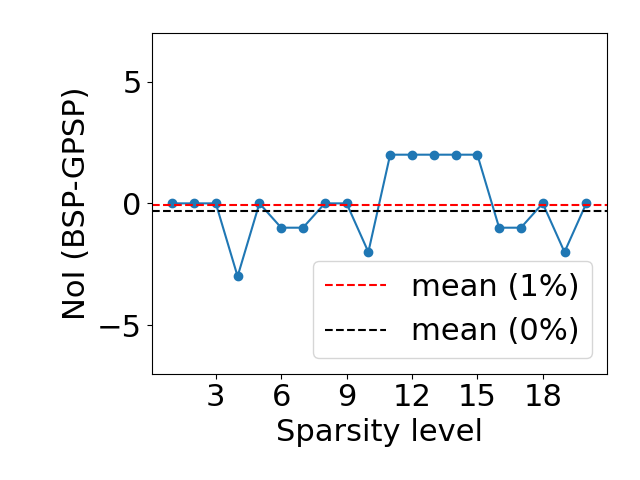}&
			\includegraphics[width=0.3\textwidth]{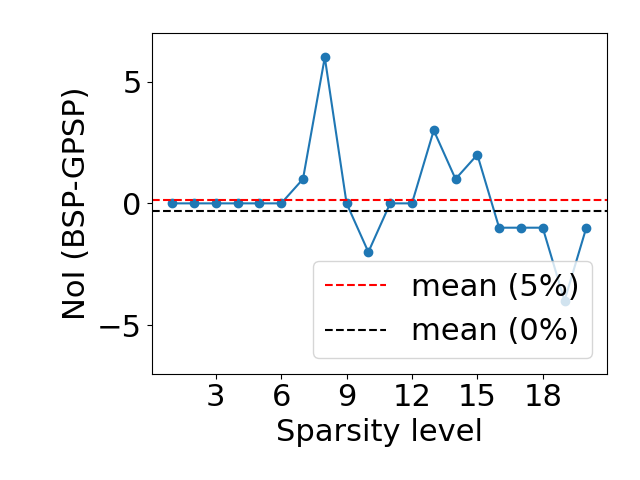}
   &
			\includegraphics[width=0.3\textwidth]{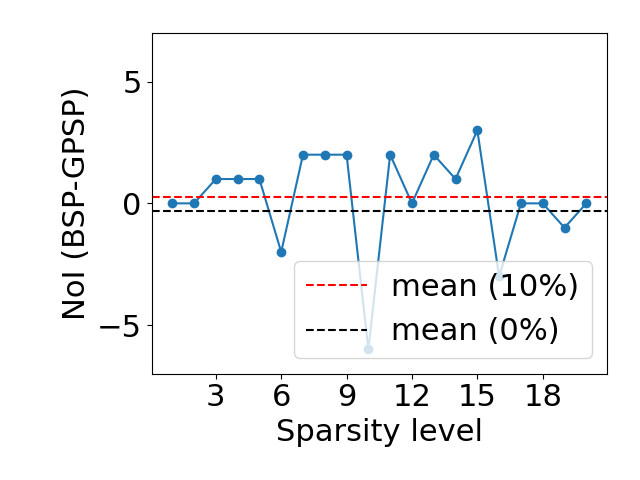}
		\end{tabular}
	\caption{For the advection-diffusion equation~\eqref{eq_exp0}, difference between the number of iterations (NoI) taken till the termination of BSP and GPSP. The blue curve represents the number of iterations of BSP minus that of GPSP. The first row shows using Dictionary I with (a) $1\%$, (b) $5\%$ and (c) $10\%$ noise. The second row with Dictionary II with (d) $1\%$, (e) $5\%$ and (f) $10\%$ noise.  The third row with Dictionary III with (g) $1\%$, (h) $5\%$ and (i) $10\%$ noise. 
Each figure shows the averaged difference between the number of BSP iterations and the number of GPSP iterations as the sparsity level varies. The dashed black line is the mean of NoI when the data is clean, and the dashed red line the mean of NoI when there is noise.
 }\label{fig_BSPGPSP_comp}
	\end{figure}
 
\section{Computational efficiency comparison between BSP and GPSP}  \label{Asec_bsp_Gpsp}

In Table~\ref{tab_eff}, we compared  the identification time for the advection-diffusion equation in~\eqref{eq_exp0} with clean data.  When $K_{\max}=10$, both BSP-IDENT and GP-IDENT are faster than the other methods, and when $K_{\max}=15$, they require more time as more candidates are generated.  We note that typically a single iteration of BSP is faster than a single iteration of GPSP, as BSP only computes vectors' inner products while GPSP involves least-square regressions. However, the speed also  depends on the number of iterations, the data, equation, and the dictionary. 
In Figure~\ref{fig_BSPGPSP_comp}, we report the difference of the number of iterations for BSP and GPSP with different noise levels and dictionaries. We observe that in general, GPSP requires fewer number of iterations than BSP when the noise level is high, and the dictionary size has an effect on this difference.

\end{document}